\documentclass[11pt]{amsart}
\usepackage{amsmath,amsfonts,amsthm,mathrsfs}
\usepackage{amssymb}

\usepackage[unicode,bookmarks]{hyperref}
\usepackage[dvipsnames]{xcolor}
\hypersetup{colorlinks=true,citecolor=NavyBlue,linkcolor=Brown,urlcolor=Orange}

\usepackage{tikz-cd,adjustbox}
\usepackage[arrow,curve,matrix]{xy}

\usepackage[margin=1in]{geometry}

\newif\ifdraft
\draftfalse

\counterwithin{equation}{subsection}
\counterwithout{subsection}{section}
\counterwithin{figure}{subsection}

\newtheorem{theorem}[equation]{Theorem}
\newtheorem*{theorem*}{Theorem}
\newtheorem{lemma}[equation]{Lemma}
\newtheorem*{lemma*}{Lemma}
\newtheorem{corollary}[equation]{Corollary}
\newtheorem{proposition}[equation]{Proposition}
\newtheorem*{proposition*}{Proposition}
\newtheorem{conjecture}[equation]{Conjecture}

\theoremstyle{definition}
\newtheorem{definition}[equation]{Definition}
\newtheorem*{definition*}{Definition}
\newtheorem{remark}[equation]{Remark}
\newtheorem*{remark*}{Remark}

\newtheorem{example}[equation]{Example}
\newtheorem*{example*}{Example}
\newtheorem*{problem*}{Problem}

\theoremstyle{plain}

\newcounter{intro}

\theoremstyle{plain}
\newtheorem{intro-conjecture}[intro]{Conjecture}
\newtheorem{intro-corollary}[intro]{Corollary}
\newtheorem{intro-theorem}[intro]{Theorem}
\newtheorem{intro-proposition}[intro]{Proposition}
\theoremstyle{definition}
\newtheorem{intro-definition}[intro]{Definition}
\newtheorem{intro-example}[intro]{Example}
\newtheorem{intro-remark}[intro]{Remark}
\newtheorem{intro-question}[intro]{Question}

\theoremstyle{definition}
\newtheorem*{unnumbered-definition}{Definition}

\newcommand{\theoremref}[1]{\hyperref[#1]{Theorem~\ref*{#1}}}
\newcommand{\lemmaref}[1]{\hyperref[#1]{Lemma~\ref*{#1}}}
\newcommand{\definitionref}[1]{\hyperref[#1]{Definition~\ref*{#1}}}
\newcommand{\propositionref}[1]{\hyperref[#1]{Proposition~\ref*{#1}}}
\newcommand{\conjectureref}[1]{\hyperref[#1]{Conjecture~\ref*{#1}}}
\newcommand{\corollaryref}[1]{\hyperref[#1]{Corollary~\ref*{#1}}}
\newcommand{\exampleref}[1]{\hyperref[#1]{Example~\ref*{#1}}}

\DeclareMathOperator{\coker}{coker}

\DeclareMathOperator{\gr}{gr}

\DeclareMathOperator{\dr}{DR}

\DeclareMathOperator{\Supp}{Supp}

\DeclareMathOperator{\codim}{codim}

\makeatletter
\let\old@caption\caption
\renewcommand*{\caption}[1]{%
	\setcounter{figure}{\value{equation}}%
	\stepcounter{equation}%
	\old@caption{#1}\relax%
}
\makeatother

\begin{document}

\title{Zero-Cycles and Holomorphic Forms in Families}

\author{Andrew Burke}
\address{Department of Mathematics, Harvard University, 1 Oxford Street, Cambridge, MA 02138, USA}
\email{aburke@math.harvard.edu}

\date{\today}

\subjclass[2010]{14B05, 14C15, 14C30, 14D06}

\begin{abstract}
We establish parallel Chow and Hodge-theoretic rigidity theorems for dominant projective morphisms between smooth complex varieties. For instance, if the general fiber satisfies $CH_0(F) = \mathbb{Z}$, then the same is true of every fiber. Similarly, if the general fiber is connected with $\mathbb{H}^0( F, \underline{\Omega}_F^p) = 0$ for $p > 0$, then the same is true of every fiber. This extends the case of rationally connected fibrations from \cite{burke dual}. We further introduce higher Chow-theoretic notions of singularities.
\end{abstract}

\maketitle

\makeatletter
\newcommand\@dotsep{4.5}
\def\@tocline#1#2#3#4#5#6#7{\relax
  \ifnum #1>\c@tocdepth 
  \else
    \par \addpenalty\@secpenalty\addvspace{#2}%
    \begingroup \hyphenpenalty\@M
    \@ifempty{#4}{%
      \@tempdima\csname r@tocindent\number#1\endcsname\relax
    }{%
      \@tempdima#4\relax
    }%
    \parindent\z@ \leftskip#3\relax
    \advance\leftskip\@tempdima\relax
    \rightskip\@pnumwidth plus1em \parfillskip-\@pnumwidth
    #5\leavevmode\hskip-\@tempdima #6\relax
    \leaders\hbox{$\m@th
      \mkern \@dotsep mu\hbox{.}\mkern \@dotsep mu$}\hfill
    \hbox to\@pnumwidth{\@tocpagenum{#7}}\par
    \nobreak
    \endgroup
  \fi}
\def\l@section{\@tocline{1}{0pt}{1pc}{}{\bfseries}}
\def\l@subsection{\@tocline{2}{0pt}{25pt}{5pc}{}}
\makeatother
\tableofcontents

\setcounter{tocdepth}{1}

\section{Introduction}

There is a close yet mysterious relationship between the Chow group of zero-cycles and spaces of holomorphic forms for complex varieties. Let $S$ be a smooth projective surface. Mumford proved \cite{mumford} that if $H^0 (S, \Omega_S^2 ) \neq 0$, then $CH_0(S)$ is infinite-dimensional in a precise sense.

Bloch conjectured that the converse holds. Indeed, let $CH_0(S)_{\mathrm{hom}}$ denote the subgroup of degree-zero cycles. Bloch's conjecture asserts that if $H^0(S, \Omega_S^2) = 0$, then the Albanese map
\begin{align*}
    \mathrm{alb}: CH_0(S)_{\mathrm{hom}} \rightarrow \mathrm{Alb}(S) := H^0(S,\Omega_S^1)^\vee/H_1(S,\mathbb Z)
\end{align*}
is an isomorphism. This is known for surfaces that are not of general type \cite{bloch conj}. A proof of the general case has recently been announced in \cite{oai}.

The conjectural Bloch--Beilinson filtration \cite{beilinson} provides a systematic framework relating Chow groups and Hodge theory; see also \cite{jannsen,voisin ii}. It predicts the following higher-dimensional form of Bloch's conjecture.

\begin{intro-conjecture}[Generalized Bloch Conjecture]\label{intro-conj BB}
Let $X$ be a smooth, connected, projective variety. Then the following conditions are equivalent:
\begin{enumerate}
    \item the degree map is an isomorphism $CH_0(X) \simeq\mathbb{Z}$; and
    \item $H^0(X, \Omega_X^p) = 0$ for every $p > 0$.
\end{enumerate}
\end{intro-conjecture}

The implication $(1)\Rightarrow(2)$ was proved by Roitman \cite{roitman}, extending Mumford's theorem to varieties of arbitrary dimension. Bloch and Srinivas later gave a different proof using their decomposition of the diagonal method \cite{bloch srinivas}. The converse remains open.

The goal of this paper is to pursue the analogy between zero-cycles and holomorphic forms by establishing parallel rigidity theorems for families. To set the stage, recall that every smooth variety admits a \emph{maximal rationally connected (MRC) fibration}, unique up to birational equivalence. Its general fiber is rationally connected, and every other rationally connected fibration factors through it. In recent work \cite{burke dual}, we established a strong rigidity theorem for rationally connected fibrations; in particular, every fiber is rationally chain connected.

Motivated by the MRC fibration, \cite{chow constant} constructs a \emph{maximal Chow trivial fibration} whose general fiber $F$ satisfies 
\[
CH_0(F) \simeq \mathbb{Z},
\]
as well as a \emph{maximal cohomologically trivial fibration} whose general fiber satisfies
\[
H^0(F,\Omega_F^p)=0
\qquad\text{for every }p>0.
\]
Conjecture \ref{intro-conj BB} predicts that these two fibrations are birationally equivalent, but this remains open. In the spirit of \cite{burke dual}, we establish rigidity theorems for both types of fibrations.

We begin with our main result on Chow trivial fibrations. Throughout, all varieties are assumed to be quasi-projective over the complex numbers. Recall that a variety $X$ is of \emph{klt type} if it is normal and supports an effective $\mathbb{Q}$-divisor $\Delta$ such that the pair $(X, \Delta)$ is klt. 

\begin{intro-theorem}\label{intro-thm CH_0}
Let $f: Y \rightarrow X$ be a dominant projective morphism, where $Y$ is irreducible and $X$ is of klt type. Suppose that the very general fiber $F$ satisfies $CH_0(F) \simeq \mathbb{Z}$. Then for any closed subvariety $Z \subset X$, there is an isomorphism
\begin{align*}
    CH_0 \big( f^{-1}(Z) \big) \simeq CH_0 \big( Z \big).
\end{align*}
In particular, every fiber $F$ satisfies $CH_0(F) \simeq \mathbb{Z}$.
\end{intro-theorem}

Applying the theorem with $Z = \{ x \}$ shows that every fiber satisfies $CH_0(F) \simeq \mathbb{Z}$, as stated above. At the other extreme, taking $Z = X$ yields an isomorphism $CH_0(Y) \simeq CH_0(X)$. When $X$ and $Y$ are also smooth, this isomorphism is known by \cite[Theorem~2.5]{chow constant}. When $f$ is a resolution of a klt singularity, it follows from \cite[Corollary $1.7$]{hm rat conn}.

Theorem \ref{intro-thm CH_0}, as well as the subsequent rigidity statements in Theorem \ref{intro-thm Omega^p}, Theorem \ref{intro-thm Omega^0}, and Conjecture \ref{intro-conj CH_0}, are known when the general fiber is rationally connected and both $X$ and $Y$ are of klt type; see \cite[Corollary $10.5$]{burke dual}. This motivated our search for stronger statements in which the general fiber is not necessarily rationally connected but nevertheless satisfies $CH_0(F) \simeq \mathbb{Z}$. 

The proof of Theorem \ref{intro-thm CH_0} has three main ingredients. First, a consequence of Tsen's theorem, recorded in Lemma \ref{lem tsen}, reduces the problem to the case $Z = \{ x \}$. We then apply the Raynaud-Gruson flatification theorem, recalled in Theorem \ref{thm flatification}, to reduce to the case in which $f$ is flat. Finally, Theorem \ref{thm flat} treats the flat case using Fulton's specialization homomorphism, recalled in Theorem \ref{thm fulton}.

Let us emphasize that Theorem \ref{intro-thm CH_0} is specific to complex varieties; see Examples \ref{exmp enriques conic} and \ref{exmp klt}. A powerful approach to rationality questions is to study the specialization of universal $CH_0$-triviality in families \cite{voisin degen, cp, totaro}. Such arguments require controlling Chow groups after base change to non-algebraically closed fields, where our methods do not apply.

To formulate a parallel Hodge-theoretic rigidity theorem, we require an appropriate notion of holomorphic forms on singular varieties. Since the sheaf of Kähler differentials can behave poorly in the singular setting, we instead use the $p$th piece 
\begin{align*}
    \mathbb{H}^0(X, \underline{\Omega}_X^p) \simeq \gr^p_F H^p(X, \mathbb{C})
\end{align*}
of the mixed Hodge filtration on cohomology. 

This space has a concrete description in terms of holomorphic forms on a resolution. Indeed, let $X' \rightarrow X$ be a resolution of singularities and let $\phi: X'' \rightarrow X' \times_X X'$ be a resolution of the fiber product. Then
\begin{align*}
    \mathbb{H}^0(X, \underline{\Omega}_X^p) \simeq \big\{ \alpha \in H^0(X' , \Omega^p_{X'}) : (\mathrm{pr}_1 \circ \phi)^* \alpha = (\mathrm{pr}_2 \circ \phi)^* \alpha \big\}.
\end{align*}
See \cite{hj} for details on differential forms in the $h$-topology. 

By extending the decomposition of the diagonal argument to singular varieties, we prove the following comparison.

\begin{intro-theorem}\label{intro-thm BB}
Let $X$ be a connected projective variety of pure dimension. If $CH_0(X) \simeq \mathbb{Z}$, then $\mathbb{H}^0(X, \underline{\Omega}_X^p) = 0$ for every $p > 0$.
\end{intro-theorem}

We expect the converse to hold as well; see Corollary \ref{cor CH_0 = Z} for supporting motivic evidence. This gives an extension of Conjecture \ref{intro-conj BB} to singular varieties. 

\begin{intro-conjecture}\label{intro-conj singular bb}
Let $X$ be a connected projective variety of pure dimension. Then the following conditions are equivalent:
\begin{enumerate}
    \item the degree map is an isomorphism $CH_0(X) \simeq\mathbb{Z}$; and
    \item $\mathbb{H}^0(X, \underline{\Omega}_X^p) = 0$ for every $p > 0$.
\end{enumerate}
\end{intro-conjecture}

We now state the Hodge-theoretic counterpart to Theorem~\ref{intro-thm CH_0}. Recall that a variety of klt type has rational singularities by a well-known result of Elkik.

\begin{intro-theorem}\label{intro-thm Omega^p}
Let $f: Y \rightarrow X$ be a dominant projective morphism, where $Y$ is irreducible and $X$ has rational singularities. Suppose that the general fiber $F$ is connected and satisfies $\mathbb{H}^0(F, \underline{\Omega}_F^p) = 0$ for every $p > 0$. Then for any closed subvariety $Z \subset X$ and $p \geq 0$, there is an isomorphism
\begin{align*}
    \mathbb{H}^0 \big( Z, \underline{\Omega}_Z^p \big) \simeq \mathbb{H}^0 \big( f^{-1}(Z), \underline{\Omega}_{f^{-1}(Z)}^p \big).
\end{align*}
In particular, every fiber $F$ satisfies $\mathbb{H}^0(F, \underline{\Omega}_F^p) = 0$ for $p > 0$.
\end{intro-theorem}

See Theorem \ref{thm Hodge singular} for a sheaf-theoretic refinement.

Both Theorem \ref{intro-thm CH_0} and Theorem \ref{intro-thm Omega^p} allow the total space $Y$ to have arbitrary singularities. Under additional regularity assumptions on $Y$, however, Theorem \ref{intro-thm Omega^p} can be strengthened. To explain this, recall that
\begin{align*}
    \mathbb{H}^0(X, \underline{\Omega}_X^p) \simeq \gr^p_F H^p(X, \mathbb{C}) \simeq \gr_F^0 \gr_p^W H^p(X, \mathbb{C})
\end{align*}
in terms of the mixed Hodge structure on cohomology. We now state a rigidity theorem for the entire Hodge graded piece
\begin{align*}
    H^p(X, \underline{\Omega}_X^0) \simeq \gr_F^0 H^p(X, \mathbb{C}),
\end{align*}
which also captures contributions from lower weights.

\begin{intro-theorem}\label{intro-thm Omega^0}
Let $f: Y \rightarrow X$ be a dominant projective morphism between connected varieties with rational singularities. Suppose that the general fiber $F$ is connected and satisfies $\mathbb{H}^p(F, \underline{\Omega}_F^0) = 0$ for every $p > 0$. Then for any closed subvariety $Z \subset X$ and $p \geq 0$, there is an isomorphism
\begin{align*}
     \mathbb{H}^p \big( Z, \underline{\Omega}_Z^0 \big) \simeq \mathbb{H}^p \big( f^{-1}(Z), \underline{\Omega}_{f^{-1}(Z)}^0 \big).
\end{align*}
In particular, every fiber $F$ satisfies $\mathbb{H}^p(F, \underline{\Omega}_F^0) = 0$ for $p > 0$.
\end{intro-theorem}

See Theorem \ref{thm Hodge general} for a sheaf-theoretic refinement which also applies when the vanishing $\mathbb{H}^p( F, \underline{\Omega}_F^0 ) = 0$ is assumed only for $p$ in a specified range. Example \ref{exmp nodal} shows that the conclusion fails when $Y$ does not have rational singularities.

Theorem \ref{intro-thm Omega^0} has in interesting consequence even when $f: Y \rightarrow X$ is a resolution of a variety with rational singularities. In this case, there is a quasi-isomorphism $Rf_* \mathscr{O}_Y \simeq \mathscr{O}_X$. A naive application of base change might therefore suggest that the fibers satisfy $H^p(F, \mathscr{O}_F) = 0$ for $p > 0$. This is false, as shown by Example \ref{exmp ratl sing}. Additionally, Example \ref{exmp reduced} shows that the vanishing can fail even after replacing $F$ with its reduction. Nevertheless, Theorem \ref{intro-thm Omega^0} proves the vanishing $\mathbb H^p(F,\underline{\Omega}_F^0)=0$ for the Du Bois complex. See Proposition \ref{prop ratl sing} for a simpler argument in this setting, which is already implicit in \cite[Lemma $2.5$]{ms}.

The proof of Theorem \ref{intro-thm Omega^0} propagates vanishing from the general fiber to every fiber by Kollár's torsion-freeness theorem for higher direct images of the canonical bundle \cite{dualizing}. It combines this torsion-freeness with general behavior of mixed Hodge modules under restriction, established in Proposition \ref{prop MHM}. The proof of Theorem \ref{intro-thm Omega^p} is similar in spirit, but instead exploits torsion-freeness of the sheaves $\mathcal{H}^0 \underline{\Omega}_X^p$ and a different restriction result, Proposition \ref{prop mhm 2}.

It is natural to ask whether Theorem \ref{intro-thm Omega^0} has a Chow-theoretic counterpart. This requires a refinement of the usual Chow group of zero-cycles for singular varieties, just as $\mathbb{H}^p(X, \underline{\Omega}_X^0)$ contains more information about lower weights of cohomology than $\mathbb{H}^0(X, \underline{\Omega}_X^p)$. Our proposed refinement is the complex of abelian groups $RCH_0(X)$ introduced in \cite{burke dual}. If $\epsilon_\bullet : X_\bullet \rightarrow X$ is a hyperresolution, then 
\begin{align*}
    RCH_0(X) \simeq \bigg( \ldots \rightarrow CH_0(X_2) \rightarrow CH_0(X_1) \rightarrow CH_0(X_0) \bigg),
\end{align*}
where the differentials are given by the alternating sum of the face maps. Its zeroth homology recovers the usual Chow group $H_0( RCH_0(X) ) \simeq CH_0(X)$. 

This leads to the following conjecture. Since $RCH_0(X)$ is difficult to study directly, it remains open even when $X$ and $Y$ are smooth.

\begin{intro-conjecture}\label{intro-conj CH_0}
Let $f: Y \rightarrow X$ be a dominant projective morphism between connected varieties of klt type. Suppose that the general fiber satisfies $CH_0(F) \simeq \mathbb{Z}$ or, equivalently, $CH_0(F)_\mathbb{Q} \simeq \mathbb{Q}$. Then for any closed subvariety $Z \subset X$, there is a quasi-isomorphism
\begin{align*}
    R CH_0 \big( f^{-1}(Z) \big)_\mathbb{Q} \simeq R CH_0 \big(Z \big)_\mathbb{Q}.
\end{align*}
In particular, every fiber $F$ satisfies $RCH_0(F)_\mathbb{Q} \simeq \mathbb{Q}[0]$, concentrated in degree zero.
\end{intro-conjecture}

The passage to rational coefficients is essential; see Example \ref{exmp enriques}. As mentioned earlier, \cite{burke dual} proves Conjecture \ref{intro-conj CH_0} when the general fiber of $f$ is rationally connected.

The interplay between the rigidity theorems above is also illuminating from a singularity-theoretic perspective. For instance, Theorem \ref{intro-thm Omega^0} applies whenever $X$ and $Y$ satisfy condition $D_0$, a slight Hodge-theoretic weakening of rational singularities. Likewise, Theorem \ref{intro-thm Omega^p} requires only that $X$ satisfy the even weaker conclusion of \cite[Theorem $0.4$]{hiatt}. Meanwhile, Theorem \ref{intro-thm CH_0} holds whenever $X$ has $CH_0$-trivial singularities, a condition arising naturally in its proof. We expect these latter two conditions to be equivalent; see \S \ref{subsection Chow Dm} for a detailed discussion.

This raises the broader question of whether one can develop a systematic Chow-theoretic hierarchy of singularities, parallel to the higher singularities in \cite{mopw, jksy, fl, mp}. As a first step, we introduce in Definition \ref{defn RCH_m} the notion of $RCH_m$-trivial singularities modeled on the Hodge-theoretic condition $D_m$ introduced in \cite{parkpopa} and, under a different name, in \cite{dor}. For $m > 0$, however, the picture becomes considerably more subtle.

\noindent
\emph{AI Disclosure.}
An AI model suggested the idea underlying Examples \ref{exmp ratl sing} and \ref{exmp reduced}, as well as Proposition \ref{prop pos ratl}. All mathematical arguments and computations were verified by the author, who assumes full responsibility for their correctness.

\noindent
{\bf Acknowledgements.}
I would like to thank Iacopo Brivio, Mihnea Popa, Stefan Schreieder, and Claire Voisin for helpful discussions on the contents of this paper.

\section{Preliminaries}

\subsection{Zero-cycles and specialization}

In this section, we recall the basic definitions and properties of Chow groups that we will need. See \cite{fulton} for further details. 

Let $X$ be a complex variety and let $i \geq 0$. Denote by $Z_i(X)$ the free abelian group generated by the classes $[V]$ of irreducible $i$-dimensional closed subvarieties $V \subset X$. If $W \subset X$ is an irreducible $(i+1)$-dimensional closed subvariety and $f \in \mathbb{C}(W)^\times$ is a rational function, then the divisor of $f$ defines a cycle class
\begin{align*}
    [ \mathrm{div}(f)] = \sum \mathrm{ord}_V(f) \cdot [V] \in Z_i(X),
\end{align*}
where the sum ranges over the $i$-dimensional subvarieties $V \subset W$. Let $\mathrm{Rat}_k(X) \subset Z_k(X)$ be the subgroup generated by these divisor classes. 

The \emph{Chow group of $i$-dimensional cycles} on $X$ is the quotient
\begin{align*}
    CH_i(X) := Z_i(X)/ \mathrm{Rat}_i(X).
\end{align*}
When $X$ is pure-dimensional, we also use the codimension grading $CH^i(X) := CH_{\dim X - i}(X)$.

Our primary object of study is the Chow group $CH_0(X)$ of zero-cycles. It is generated by the classes of closed points of $X$, modulo the divisor classes of rational functions on irreducible curves. If $X$ is projective, there is a surjective degree map 
\begin{align*}
    \mathrm{deg}: CH_0(X) \rightarrow \mathbb{Z}, \qquad \sum_i n_i[x_i] \mapsto \sum_i n_i.
\end{align*}
We will often write $CH_0(X) \simeq \mathbb{Z}$ to indicate that this degree map is an isomorphism. 

We record the following observation for future use. It is the key reason why Theorem \ref{intro-thm CH_0} requires no regularity assumptions on the total space; see Remark \ref{rmk Y irred}.

\begin{lemma}\label{lem dominant}
Let $f: Y \rightarrow X$ be a dominant projective morphism. Then the induced pushforward map 
\begin{align*}
    f_*: CH_0(Y) \rightarrow CH_0(X) 
\end{align*}
is surjective. In particular, if $CH_0(Y) \simeq \mathbb{Z}$, then $CH_0(X) \simeq \mathbb{Z}$.
\end{lemma}

\begin{proof}
Since $f$ is projective and dominant, it is surjective. The group $CH_0(X)$ is generated by the classes of closed points. Given a closed point $x \in X$, choose a closed point $y \in Y$ with $f(y) = x$. Then $f_*[y] = [x]$.
\end{proof}

We will consider morphisms whose fibers satisfy $CH_0(X) \simeq \mathbb{Z}$. The following lemma shows that this condition can be checked after tensoring with $\mathbb{Q}$.

\begin{lemma}\label{lem CH_0 Q}
Let $X$ be a complex projective variety. If $CH_0(X)_\mathbb{Q} \simeq \mathbb{Q}$, then $CH_0(X) \simeq \mathbb{Z}$.
\end{lemma}

\begin{proof}
Fix a closed point $x_0 \in X$. The degree map $\deg: CH_0(X) \rightarrow \mathbb{Z}$ is split by the assignment $n \mapsto n [x_0]$, so there is a direct sum decomposition
\begin{align*}
    CH_0(X) \simeq CH_0(X)_{\mathrm{hom}} \oplus \mathbb{Z}.
\end{align*}
The hypothesis $CH_0(X)_\mathbb{Q} \simeq \mathbb{Q}$ therefore implies $CH_0(X)_{\mathrm{hom}, \mathbb{Q}} = 0$. Equivalently, every element of $CH_0(X)_{\mathrm{hom}}$ is torsion.

For each $m \geq 1$, define
\begin{align*}
    V_m := \bigg \{ x \in X : m \big( [x] - [x_0] \big) = 0 \in CH_0(X) \bigg \}.
\end{align*}
Since $CH_0(X)_{\mathrm{hom}}$ is torsion, we have $X = \bigcup_{m \geq 1} V_m$.

We claim that each $V_m$ is a countable union of Zariski closed subsets. Indeed, if $\Delta_m: X \rightarrow \mathrm{Sym}^m X \times \mathrm{Sym}^m X$ denotes the map $x \mapsto (mx, m x_0)$, then $V_m = \Delta_m^{-1}( W_m)$, where
\begin{align*}
    W_m = \bigg \{ (y,z) \in \mathrm{Sym}^m X \times \mathrm{Sym}^m X : [y] = [z] \in CH_0(X) \bigg \}.
\end{align*}
By \cite[Lemma $10.7$]{voisin ii}, $W_m$ is a countable union of Zariski closed subsets, and hence so is $V_m$.

Let $X_i$ be an irreducible component of $X$. It cannot be expressed as a countable union of proper closed subvarieties, so $X_i \subset V_{m_i}$ for some $m_i \geq 1$. Taking $m$ to be a common multiple of the $m_i$, we obtain $X = V_m$. Since $CH_0(X)_{\mathrm{hom}}$ is generated by classes of the form $([x] - [x_0])$, this gives
\begin{align*}
    m \cdot CH_0(X)_{\mathrm{hom}} = 0.
\end{align*}

But we claim that $CH_0(X)_{\mathrm{hom}}$ is also a divisible group. Assuming this for the moment, we have $CH_0(X)_{\mathrm{hom}} = m \cdot CH_0(X)_{\mathrm{hom}} = 0$. Hence, $CH_0(X) \simeq \mathbb{Z}$, completing the proof.

To prove this claim, note that $CH_0(X)_{\mathrm{hom}}$ is generated by classes $([x] - [y]),$ where $x$ and $y$ lie in the same irreducible component $X_i$ in $X$. For any such class, choose an irreducible curve $C \subset X$ passing through $x$ and $y$. Let $\tilde{x}$ and $\tilde{y}$ be lifts to the normalization $\tilde{C}$ of $C$. Then the class $([\tilde{x}] - [\tilde{y}]) \in CH_0(\tilde{C})_{\mathrm{hom}} \simeq \mathrm{Jac}(\tilde{C})$ is divisible, hence so is $[x] - [y]$.
\end{proof}

Next, we recall the specialization homomorphism for Chow groups, which is our main tool for analyzing zero-cycles in flat families in the proof of Theorem \ref{thm flat}.

\begin{theorem}\label{thm fulton}\cite[\S 4.4]{fulton specialization}
Let $R$ be a discrete valuation ring with fraction field $K$ and residue field $k$. Let $X$ be a flat quasi-projective $R$-scheme, and denote its generic and special fibers by $X_K = X \otimes_R K$ and $X_k = X \otimes_R k$, respectively. 

Then, for every $i \geq 0$, there exists a specialization homomorphism
\begin{align*}
    \sigma_X: CH_i(X_K) \rightarrow CH_i(X_k)
\end{align*}
defined as follows. Let $V \subset X_K$ be an irreducible closed subvariety of dimension $i$, and let $\overline{V} \subset X$ denote its closure. Then $\sigma_X([V]) = [\overline{V}_k]$, where $\overline{V}_k = \overline{V} \times_{ \mathop{\mathrm{Spec}} R } \mathop{\mathrm{Spec}} k$ and $[\overline{V}_k]$ denotes the associated $i$-dimensional cycle, including its scheme-theoretic multiplicities.
\end{theorem}

We will also use the complex of abelian groups $RCH_0(X)$ introduced in \cite{burke dual}. Given a hyperresolution $\epsilon_\bullet: X_\bullet \rightarrow X$, there is a quasi-isomorphism
\begin{align*}
    RCH_0(X) \simeq \bigg( \ldots \rightarrow CH_0(X_2) \rightarrow CH_0(X_1) \rightarrow CH_0(X_0) \bigg).
\end{align*}
As with the usual Chow group of zero-cycles, $RCH_0(X)$ is covariant for proper morphisms. Its zeroth homology recovers the usual Chow group:

\begin{lemma}\cite[Proposition $6.13$]{burke dual}\label{lem RCH_0}
For every variety $X$, there is an isomorphism
\begin{align*}
    H_0 \big( RCH_0(X) \big) \simeq CH_0(X).
\end{align*}
\end{lemma}

\subsection{Mixed Hodge Modules and de Rham complex}

In this section, we review the theory of mixed Hodge modules, introduced by Saito in \cite{saito 88, saito 90}.

Let $X$ be a complex variety embedded in a smooth variety $V$. Then a \emph{mixed Hodge module} $M \in \mathrm{MHM}(X)$ is the data of a tuple
\begin{align*}
    \big( \mathcal{M}, F_\bullet, W_\bullet, (K, W_\bullet), \alpha \big).
\end{align*}
Here, $\mathcal{M}$ is a regular holonomic $\mathcal{D}_V$-module supported on $X$, $F_\bullet =F _\bullet \mathcal{M}$ is a good filtration, and $W_\bullet = W_\bullet \mathcal{M}$ is a finite filtration by $\mathcal{D}_V$-modules. The pair $(K, W_\bullet)$ consists of an algebraically constructible $\mathbb{Q}$-perverse sheaf on $X^{\mathrm{an}}$ with a finite increasing filtration $W_\bullet = W_\bullet K$, and $\alpha$ is a filtered isomorphism 
\begin{align*}
    \alpha:\mathbb{C} \otimes_\mathbb{Q} (K, W) \rightarrow \mathrm{DR}_X(\mathcal{M}, W)
\end{align*}
We make the convention that all $\mathcal{D}$-modules act on the left. However, we will index the Hodge filtration following the conventions for right $\mathcal{D}$-modules. 

The data defining a mixed Hodge module is required to satisfy many compatibility conditions, which we will not review here. The idea is that a mixed Hodge module on a point is precisely a graded-polarizable mixed Hodge structure, with Hodge filtration $F_\bullet$ and weight filtration $W_\bullet$. On higher-dimensional varieties, the definition proceeds inductively on dimension, imposing compatibility conditions through nearby and vanishing cycles associated to local functions.

The category $\mathrm{MHM}(X)$ of mixed Hodge modules on $X$ is abelian, and we denote its bounded derived category by $D^b(\mathrm{MHM}(X))$. A major accomplishment of Saito's theory is to endow this category with a six functor formalism.

\begin{theorem}\cite[Theorem $0.1$]{saito 90}
If $X$ is a complex variety, there is an exact dualizing functor
\begin{align*}
    \mathbb{D}: \mathrm{MHM}(X) \rightarrow \mathrm{MHM}(X)^{\mathrm{op}}.
\end{align*}
Moreover, if $f: Y \rightarrow X$ is a morphism of complex varieties, there exist natural functors
\begin{align*}
    f_*, f_!: D^b(\mathrm{MHM}(Y)) \rightarrow D^b(\mathrm{MHM}(X)), \hspace{.3cm}  f^*, f^!: D^b(\mathrm{MHM}(X)) \rightarrow D^b(\mathrm{MHM}(Y))
\end{align*}
with $f^*$ left adjoint to $f_*$ and $f_!$ left adjoint to $f^!$. These are compatible with the corresponding functors on complexes of constructible sheaves.
\end{theorem}

The functors satisfy the usual compatibility isomorphisms. For example, $f_* \circ \mathbb{D} \simeq \mathbb{D} \circ f_!$ and $f^* \circ \mathbb{D} \simeq \mathbb{D} \circ f^!$. We will generally use these without further comment. Let us record, however, the following base-change isomorphisms.

\begin{proposition}\cite[\S $4$]{saito 90}\label{prop base change}
Consider a cartesian square
\[
\begin{tikzcd}
Y' \arrow[r,"g'"] \arrow[d,"f'"']
  & Y \arrow[d,"f"] \\
X' \arrow[r,"g"']
  & X.
\end{tikzcd}
\]
Then there are natural isomorphisms $g^! \circ f_* \simeq f'_* \circ (g')^!$ and $g^* \circ f_! \simeq (f')_! \circ (g')^*$.
\end{proposition}

We also need the localization triangles associated with complementary closed and open embeddings.

\begin{proposition}\cite[$4.4.1$]{saito 90}\label{prop dist tri}
Let $i: Z \hookrightarrow X$ be a closed embedding, and let $j: U = X \setminus Z \rightarrow X$ be its complementary open embedding. Then for any $M^\bullet \in D^b(\mathrm{MHM}(X))$, there are distinguished triangles
\begin{align*}
    j_! j^! M^\bullet \rightarrow M^\bullet \rightarrow i_* i^* M^\bullet \xrightarrow[]{+1}
\end{align*}
and
\begin{align*}
    i_! i^! M^\bullet \rightarrow M^\bullet \rightarrow j_* j^* M^\bullet \xrightarrow[]{+1}
\end{align*}
Moreover, $i_* = i_!$ and $j^* = j^!$, and these two functors are exact.
\end{proposition}

Suppose that $X$ is embedded in a smooth variety $V$ of dimension $d$, and let $M \in \mathrm{MHM}(X)$. Denote by $\mathcal{M}$ the underlying $\mathcal{D}_{V}$-module. The \emph{de Rham complex} of $M$ is
\begin{align*}
    \mathrm{DR}_X(\mathcal{M}) := \bigg( \mathcal{M} \xrightarrow{\nabla} \Omega_{V}^1 \otimes_\mathscr{O} \mathcal{M} \xrightarrow{\nabla} \ldots \xrightarrow{\nabla} \Omega_{V}^d \otimes_\mathscr{O} \mathcal{M} \bigg).
\end{align*}
The Hodge filtration $F_\bullet \mathcal{M}$ induces a filtration on $\mathrm{DR}_X(\mathcal{M})$. Its $p$th associated graded piece is of the form
\begin{align}\label{eqn grDR}
    \gr_{p}^F \dr(M) := \bigg( \gr_{p-d}^F \mathcal{M} \rightarrow \Omega_{V}^1 \otimes_\mathscr{O} \gr_{p-d+1}^F \mathcal{M} \rightarrow \ldots \rightarrow \Omega_{V}^d \otimes_\mathscr{O} \gr_p^F \mathcal{M} \bigg).
\end{align}

This is a complex of coherent sheaves on $X$, and its quasi-isomorphism class depends only on $M$, not on the chosen embedding $X \hookrightarrow V$. More generally, this construction extends to an exact functor
\begin{align*}
    \gr_p^F \dr: D^b(\mathrm{MHM}(X)) \rightarrow D^b_{\mathrm{coh}}(\mathscr{O}_X),
\end{align*}
where $D^b_{\mathrm{coh}}(\mathscr{O}_X)$ denotes the bounded derived category of coherent sheaves on $X$. 

The graded de Rham functor is compatible with duality and proper pushforward. To formulate these compatibilities, consider the Grothendieck duality functor
\begin{align*}
    \mathbb{D}_X: D^b_{\mathrm{coh}}(\mathscr{O}_X) \rightarrow D^b_{\mathrm{coh}}(\mathscr{O}_X)^{\mathrm{op}}, \qquad \mathscr{F}^\bullet \mapsto \mathrm{R}\mathcal{H}om_{\mathscr{O}_X}(\mathscr{F}^\bullet, \omega_X^\bullet[n]).
\end{align*}
Here, $n = \dim X$ and $\omega_X^\bullet$ is the dualizing complex. With this normalization, if $X$ is smooth, then $\omega_X^\bullet [n] = \omega_X$. We have the following result.

\begin{proposition}\cite[Proposition $3.1$ and Lemma $3.2$]{park}\label{prop D comp}
Let $X$ be a complex variety with $n = \dim X$. If $M^\bullet \in D^b(\mathrm{MHM}(X))$ and $p \in \mathbb{Z}$, then
\begin{align*}
    \gr_p^F \dr \big( \mathbb{D} M^\bullet \big) \simeq \bigg( \mathbb{D}_X \gr_{-p}^F \dr \big( M^\bullet \big) \bigg)[n].
\end{align*}
Moreover, if $f: X \rightarrow Y$ is a projective morphism, then
\begin{align*}
    \gr_p^F \dr \big( f_* M^\bullet \big) \simeq Rf_* \gr_p^F \dr \big( M^\bullet \big).
\end{align*}
\end{proposition}

By contrast, the graded de Rham complex has more subtle behavior under pullback. In \S \ref{subsection restriction}, we study this behavior for a closed embedding $i: Z \rightarrow X$.

\subsection{Du Bois complex and higher singularities}
For a detailed discussion of the material presented here, see \cite[Ch. B]{parkpopa}.

Let $X$ be a complex variety of dimension $n$. For each $0 \leq p \leq n$, the \emph{$p$-th Du Bois complex} of $X$ is an object $\underline{\Omega}_X^p$ in the derived category of coherent sheaves,
introduced in \cite{du bois}. It provides a Hodge-theoretic replacement for singular varieties of the sheaf $\Omega_X^p$ of holomorphic $p$-forms.

The Du Bois complex admits an interpretation in terms of mixed Hodge modules. Indeed, if $a: X \rightarrow \mathop{\mathrm{Spec}} \mathbb{C}$ is the structure morphism, set $\mathbb{Q}_X^H := a^* \mathbb{Q}^H$. Saito \cite[Theorem $0.2$]{saito 00} proved that there is a natural isomorphism
\begin{align}\label{eqn db mhm}
    \underline{\Omega}_X^p \simeq \gr_{-p}^F \dr \big( \mathbb{Q}_X^H[n] \big)[p-n].
\end{align}
We will make this identification throughout the paper.

We next recall Deligne's mixed Hodge structure on cohomology in the case where $X$ is projective. As in the smooth setting, there is a Hodge-de Rham spectral sequence
\begin{align}\label{Hodge-dR ss}
    E_1^{p,q} = \mathbb{H}^q(X, \underline{\Omega}_X^p) \Rightarrow H^{p+q}(X, \mathbb{C}),
\end{align}
which degenerates at the $E_1$-page; see e.g. \cite[Proposition $7.1$]{peters steenbrink}. The induced filtration on $H^{p+q}(X, \mathbb{C})$ agrees with the Hodge filtration. Consequently,
\begin{align*}
    \gr_{-p}^F H^{p+q}(X, \mathbb{C}) \simeq \mathbb{H}^q(X, \underline{\Omega}_X^p).
\end{align*}

The weight filtration can be described using a hyperresolution $\epsilon_\bullet : X_\bullet \rightarrow X$. Indeed, there is a cohomological weight spectral sequence
\begin{align*}
    E_1^{p,q} = H^q(X_p, \mathbb{Q}) \Rightarrow H^{p+q}(X, \mathbb{C}).
\end{align*}
It degenerates at the $E_2$-page, and the induced filtration on cohomology is the weight filtration. In particular, its associated graded pieces are given by
\begin{align}\label{eqn weights}
    \gr_q^W H^{p+q}(X, \mathbb{Q}) = \frac{\ker \bigg( H^q(X_p, \mathbb{Q}) \rightarrow H^q(X_{p+1}, \mathbb{Q}) \bigg)}{\mathrm{im} \bigg( H^q(X_{p-1}, \mathbb{Q}) \rightarrow H^q(X_{p}, \mathbb{Q}) \bigg)}.
\end{align}
Moreover, these maps are strict with respect to the Hodge filtration. See \cite[\S 2]{bps} for further details.

We now introduce a Hodge-theoretic notion of singularities in terms of the Du Bois complex. Condition $D_m$, studied in \cite{parkpopa} and, under a different name, in \cite{dor}, stands for duality through level $m$. 

A pure $n$-dimensional variety $X$ is said to satisfy \emph{condition $D_m$} if the canonical morphism
\begin{align*}
    \underline{\Omega}_X^p \rightarrow \mathbb{D}_X(\underline{\Omega}_X^{n-p})
\end{align*}
is an isomorphism for every $0 \leq p \leq m$.

There is a concrete characterization of condition $D_m$ in terms of the Hodge filtration on local cohomology.

\begin{theorem}\cite[Theorem D(2)]{dor}\label{thm Dm}
Let $X$ be a purely $n$-dimensional variety. Then $X$ satisfies condition $D_m$ if and only if, for every point $x \in X$,
\begin{align*}
    F_{m-n} H^j_{ \{ x \} }(X, \mathbb{C}) = 
                \begin{cases}
                    \mathbb{C}, & \text{ if } j = 2n, \\
                    0, & \text{ if } j \neq 2n.
                \end{cases}
\end{align*}
\end{theorem}

We will be primarily interested in the case $m = 0$, where the condition is more concrete.

\begin{lemma}\label{lem D0}
Let $X$ be an $n$-dimensional variety with resolution of singularities $f: Y \rightarrow X$. Then there are quasi-isomorphisms
\begin{align*}
    \underline{\Omega}_X^n \simeq Rf_* \omega_{Y} \simeq f_* \omega_{Y}.
\end{align*}

Therefore, $X$ satisfies condition $D_0$ if and only if $\underline{\Omega}_X^0 \simeq Rf_* \mathscr{O}_Y$ or, dually, $\mathbb{D} \underline{\Omega}_X^0 \simeq f_* \omega_Y$.
\end{lemma}

\begin{proof}
The first isomorphism $\underline{\Omega}_X^n \simeq Rf_* \omega_{Y}$ follows from the description of the Du Bois complex in terms of hyperresolutions. Since $f$ is a resolution of singularities, Grauert-Riemenschneider vanishing gives $R^i f_* \omega_Y = 0$ for all $i > 0$, and hence $Rf_* \omega_{Y} \simeq f_* \omega_Y$.
\end{proof}

Condition $D_0$ is weaker than the classical notion of rational singularities. More precisely, a variety has rational singularities if and only if it has Du Bois singularities and satisfies condition $D_0$.

For higher values of $m$, this characterization extends to the higher notions of Du Bois and rational singularities developed in \cite{mopw, jksy, fl, mp}. Combining \cite[Proposition $9.4$]{PSV} with \cite[\S $4.5$]{dor} yields the following result.

\begin{theorem}
A normal variety $X$ has pre-$m$-rational singularities if and only if it has pre-$m$-Du Bois singularities and satisfies condition $D_m$.
\end{theorem}

There is an equivalent formulation of condition $D_m$ in terms of intersection cohomology. This requires some more terminology. Define the \emph{intersection complex Hodge module of $X$} by the formula
\begin{align*}
   \mathrm{IC}_X^H := \gr_n^W \mathcal{H}^0 \mathbb{Q}_X^H[n].
\end{align*}

Its hypercohomology endows the intersection cohomology $IH^\bullet(X, \mathbb{Q})$ with respect to the middle perversity with a pure Hodge structure. Moreover, following \cite[Definition $6.1$]{parkpopa}, the \emph{RHM-defect object of $X$} is the object $\mathcal{K}_X^\bullet \in D^b(\mathrm{MHM}(X))$ sitting in the natural distinguished triangle
\begin{align}\label{eqn K}
    \mathcal{K}_X^\bullet \rightarrow \mathbb{Q}_X^H[n] \rightarrow \mathrm{IC}_X^H \xrightarrow[]{+1}
\end{align}

In analogy with the description of the Du Bois complex in \eqref{eqn db mhm}, the \emph{intersection Du Bois complex}, introduced in \cite[Definition $3.3$]{parkpopa}, is defined by
\begin{align*}
    I \underline{\Omega}_X^p := \gr_{-p}^F \dr \big( \mathrm{IC}_X^H \big) [p-n].
\end{align*}
Then condition $D_m$ is equivalent to the natural morphism
\begin{align*}
    \underline{\Omega}_X^p \rightarrow I \underline{\Omega}_X^p
\end{align*}
being an isomorphism for every $0 \leq p \leq m$; see \cite[Remark $3.8$]{parkpopa} for details. Equivalently, $\gr_{-p}^F \dr ( \mathcal{K}_X^\bullet ) = 0$ for all $0 \leq p \leq m$.

We will need the following slight generalization of a result of Hiatt for the proof of Theorem~\ref{intro-thm Omega^p}. We include an argument using the mixed Hodge module language of this paper.

\begin{proposition}\cite[Theorem $0.4$]{hiatt}\label{prop hiatt}
If $X$ satisfies condition $D_m$, then for every $0 \leq p \leq n$, the natural map
\begin{align*}
    \mathcal{H}^i ( \underline{\Omega}_X^p) \rightarrow \mathcal{H}^i ( I \underline{\Omega}_X^p )
\end{align*}
is an isomorphism for $i \leq m$ and injective for $i = m+1$.
\end{proposition}

For lack of a better name, we refer to the conclusion of this proposition as \emph{condition $E_m$}. Thus, Proposition \ref{prop hiatt} says that condition $D_m$ implies condition $E_m$.

\begin{proof}
By assumption, we have vanishing $\gr_{-p}^F \dr ( \mathcal{K}_X^\bullet) = 0$ for all $p \leq m$. Equivalently, by the dual of \cite[Corollary $6.3$]{burke lcdef},
\begin{align*}
    \gr_{-p}^F \dr \big( \mathcal{H}^i \mathcal{K}_X^\bullet \big) = 0 \hspace{.4cm} \text{ for all } i \text{ and } p \leq m.
\end{align*}
By \cite[Proposition $6.4$]{parkpopa}, the mixed Hodge module $\mathcal{H}^i \mathcal{K}_X^\bullet$ has weights $\leq (i+n-1)$. Combining this fact with the duality result of \cite[Proposition $5.2$]{parkpopa}, we obtain vanishing
\begin{align*}
    \gr_{-p}^F \dr \big( \mathcal{H}^i \mathcal{K}_X^\bullet \big) = 0 \hspace{.4cm} \text{ whenever } i \leq p - n + m +1.
\end{align*}

Then Corollary \ref{cor cohom vs grDR} proved independently below shows that the latter vanishing is equivalent to
\begin{align*}
    \mathcal{H}^i \gr_{-p}^F \dr \big( \mathcal{K}_X^\bullet \big) = 0 \hspace{.4cm} \text{ whenever } i \leq p - n + m +1.
\end{align*}
Applying $\gr_{-p}^F \dr(-)$ to the distinguished triangle \eqref{eqn K}, we see that this is equivalent to the conclusion of the proposition.
\end{proof}

\section{Chow and Hodge Rigidity in families}

\subsection{Rigidity of zero-cycles}

The goal of this section is to study the Chow group of zero-cycles in families, culminating in a proof of Theorem \ref{intro-thm CH_0}.

We begin by isolating the condition appearing in the conclusion of Theorem \ref{intro-thm CH_0} and recording several equivalent formulations. In particular, that condition can be tested on closed fibers.

\begin{lemma}\label{lem CH_0 trivial}
Let $f: Y \rightarrow X$ be a dominant projective morphism of varieties. Then the following conditions are equivalent:
\begin{enumerate}
    \item $CH_0 \big(f^{-1}(x) \big) \simeq \mathbb{Z}$ for every closed point $x \in X$;
    \item $CH_0 \big(f^{-1}(x) \big)_\mathbb{Q} \simeq \mathbb{Q}$ for every closed point $x \in X$;
    \item $CH_0 \big(f^{-1}(Z) \big) \simeq CH_0(Z)$ for every closed subvariety $Z \subset X$; and 
    \item $CH_0 \big(f^{-1}(Z) \big)_\mathbb{Q} \simeq CH_0(Z)_\mathbb{Q}$ for every closed subvariety $Z \subset X$.
\end{enumerate}
When these conditions hold, we say that $f$ is \emph{$CH_0$-trivial}. If $f$ is a resolution of singularities, we say that $X$ has \emph{$CH_0$-trivial singularities}.
\end{lemma}

Lemma \ref{lem wft} below shows that, when $f$ is a resolution of singularities, this condition is independent of the chosen resolution. The terminology of $CH_0$-triviality was introduced in the more general framework of \cite[Definition $7.1$]{burke dual}. However, for the broader class of invariants considered there, the conditions of Lemma \ref{lem CH_0 trivial} need not be equivalent, and the term is reserved for condition $(3)$. 

The proof of Lemma \ref{lem CH_0 trivial} relies on the following consequence of Tsen's theorem.

\begin{lemma}\cite[Propositions $10, 11$]{serre}\label{lem tsen}
Let $L/K$ be a finite extension of function fields of complex curves. Then the norm map $N = N_{L/K}: L^\times \rightarrow K^\times$ is surjective.
\end{lemma}

\begin{proof}
Let $a \in K^\times$. We will find $b \in L^\times$ such that $N(b) = a$. Recall Tsen's theorem asserts that $K$ is a $C_1$ field, that is every homogeneous polynomial $F \in K[X_1, \ldots , X_n]$ of degree $d \leq n-1$ has a nontrivial zero.

Set $d = [L:K]$ and choose a $K$-basis $e_1, \ldots e_d$ of $L$. For a generic linear combination $X = X_1 e_1 + \ldots X_d e_d$, let $m_X: L \rightarrow L$ denote multiplication by $X$. Then the norm
\begin{align*}
    N(X) = \det( m_X),
\end{align*}
is a homogeneous polynomial of degree $d$ in the variables $X_1, \ldots , X_d$. Indeed, if we write $e_i \cdot e_j = \sum_{k = 1}^d a^{ij}_k e_k$ for constants $a^{ij}_k \in K$, it is easy to see that $m_X$ is given by the matrix with entries $(m_X)_{ij} = \sum_{\ell =1}^d X_\ell a_i^{\ell j}$, and these are linear in $X_1, \ldots , X_d$.

Consider now the homogeneous polynomial $F(X_1, \ldots X_d, Y) = N(X) - a Y^d$. It is homogeneous of degree $d$ in $d+1$ variables. By Tsen's theorem, it admits a nontrivial zero $(x_1, \ldots, x_d, y)$. Since $N$ is a non-degenerate form, $y \neq 0$. After rescaling, we may assume $y = 1$. Then the element $b= x_1 e_1 + \ldots x_d e_d \in L^\times$ satisfies $N(b) = a$.
\end{proof}

We are now prepared for the proof of Lemma \ref{lem CH_0 trivial}.

\begin{proof}[Proof of Lemma \ref{lem CH_0 trivial}]
By Lemma \ref{lem CH_0 Q}, conditions $(1)$ and $(2)$ are equivalent. Moreover, it is clear that $(3) \Rightarrow (4) \Rightarrow (2)$. It suffices to prove that $(1) \Rightarrow (3)$.

Assume $(1)$, and fix a closed subvariety $Z\subset X$. Define a homomorphism
\begin{align*}
    \theta: Z_0(Z) \rightarrow CH_0(f^{-1}(Z))
\end{align*}
on the group $Z_0(Z)$ of zero-cycles of $Z$ as follows. For a closed point $p \in Z$, choose any point $q \in f^{-1}(p)$, and then set $\theta([p]) = [q]$. Since $CH_0(f^{-1}(p)) \simeq \mathbb{Z}$, the class $[q]$ is independent of the choice of pre-image, so $\theta$ is well-defined. We claim that $\theta$ respects rational equivalence and hence descends to a homomorphism
\begin{align*}
    \theta: CH_0(Z) \rightarrow CH_0(f^{-1}(Z)).
\end{align*}
Once this is established, it is clear that $\theta$ is an inverse to the induced pushforward map $f_*: CH_0(f^{-1}(Z)) \rightarrow CH_0(Z)$.

To prove that $\theta$ respects rational equivalence, let $C \subset Z$ be an irreducible curve and let $\phi \in \mathbb{C}(C)^\times$ be a rational function. Choose an irreducible curve $D \subset f^{-1}(C)$ such that the restriction $f|_D: D \rightarrow C$ is finite. Passing to function fields $K = \mathbb{C}(C)$ and $L = \mathbb{C}(D)$, we obtain a finite extension $L/K$. By Lemma \ref{lem tsen}, the norm map $N=N_{L/K}$ is surjective, so there exists a rational function $\psi \in L^\times$ for which $N_{L/K}(\psi) = \phi$. The compatibility of proper pushforward with the norm map gives $f_* [ \mathrm{div}(\psi) ] = [\mathrm{div}(\phi)]$. Hence,
\begin{align*}
    \theta \big( \mathrm{div(\phi) \big) = \mathrm{div}(\psi)}
\end{align*}
lies in $\mathrm{Rat}_0(f^{-1}(Z))$. We conclude that $\theta$ respects rational equivalence, and this proves $(3)$.
\end{proof}

We obtain the following immediate consequence.

\begin{lemma}\label{lem wft}
Let $f: Y \rightarrow X$ be a projective birational morphism between smooth varieties. Then $f$ is $CH_0$-trivial.
\end{lemma}

\begin{proof}
By the weak factorization theorem \cite{weak factorization}, $f$ factors as a finite sequence of blow-ups and blow-downs along smooth centers. It therefore suffices to treat the case in which $Y = \mathrm{Bl}_C X$ is the blow-up of $X$ along a smooth center $C$. In this situation, each fiber of $f$ is either a point or a copy of projective space $\mathbb{P}^{\codim_X(C)-1}$. In either case, $CH_0(f^{-1}(x)) \simeq \mathbb{Z}$. The claim therefore follows from condition $(1)$ in Lemma \ref{lem CH_0 trivial}.
\end{proof}

We now state the main result of this section.

\begin{theorem}[{\(\!\Rightarrow\) Theorem~\ref{intro-thm CH_0}}]\label{thm CH_0}
Let $f: Y \rightarrow X$ be a dominant projective morphism, where $Y$ is irreducible and $X$ has $CH_0$-trivial singularities. Suppose that the very general fiber $F$ satisfies $CH_0(F) \simeq \mathbb{Z}$. Then for any closed subvariety $Z \subset X$, the induced map
\begin{align*}
    CH_0 \big( f^{-1}(Z) \big ) \rightarrow CH_0 \big( Z \big)
\end{align*}
is an isomorphism.
\end{theorem}

This is slightly stronger than Theorem \ref{intro-thm CH_0}. Indeed, \cite[Corollary H]{burke dual} shows that every variety of klt type has $CH_0$-trivial singularities. The $CH_0$-triviality condition is the natural hypothesis for the argument below.

\begin{remark}[{\bf Singularities of $Y$}]\label{rmk Y irred}
It is noteworthy that Theorem \ref{thm CH_0} imposes no restrictions on the singularities of $Y$. However, when the singular locus of $Y$ does not dominate $X$, it is easy to reduce to the case of a smooth total space.

Indeed, let $s: Y' \rightarrow Y$ be a resolution of singularities which is an isomorphism on the smooth locus. The very general fibers of $f \circ s$ and $f$ are isomorphic, so, assuming Theorem \ref{thm CH_0} for smooth total spaces, we may apply it to the composition $f \circ s$. Hence, for every closed subvariety $Z \subset X$, the composition
\begin{align*}
    CH_0\big( (f \circ s)^{-1}(Z) \big) \xrightarrow{s_*} CH_0 \big( f^{-1}(Z) \big) \xrightarrow{f_*} CH_0 \big( Z \big)
\end{align*}
is an isomorphism. It follows that $s_*$ is injective. On the other hand, since $s$ is a dominant morphism, $s_*$ is surjective by Lemma \ref{lem dominant}. It follows that $s_*$ is an isomorphism, and consequently so is $f_*$.

Our proof of Theorem \ref{thm CH_0} below does not require any such reduction.
\end{remark}

\begin{remark}[{\bf Very general fiber}]
The hypothesis on a very general fiber is used to show that the geometric generic fiber has $CH_0 \simeq \mathbb{Z}$. Although this is weaker than assuming $CH_0\simeq\mathbb{Z}$ for every general fiber, it is precisely the hypothesis needed for the passage to the geometric generic fiber.
\end{remark}

The strategy for Theorem \ref{thm CH_0} is to reduce to the case where $f$ flat by the following version of the Raynaud-Gruson flatification theorem, and then to invoke Theorem \ref{thm flat}.

\begin{theorem}\cite{flatification}\label{thm flatification}
Let $f: Y \rightarrow X$ be a dominant morphism of varieties, with $X$ irreducible. Then there exists a blow-up $r: X' \rightarrow X$ along a closed subscheme not containing the generic point of $X$, such that, if $Y' \subset Y \times_X X'$ denotes the strict transform of $Y$, then the induced morphism $f': Y' \rightarrow X'$ is flat.
\end{theorem}

We are now prepared to prove Theorem \ref{thm CH_0}.

\begin{proof}[Proof of Theorem \ref{thm CH_0}]
Choose a flatification of $f$ as in Theorem \ref{thm flatification}. We obtain a commutative square
\[
\begin{tikzcd}
Y' \arrow[r, "s"] \arrow[d, "f'"'] & Y \arrow[d, "f"] \\
X' \arrow[r, "r"'] & X
\end{tikzcd}
\]
where $f'$ is flat and $r$ is birational. Since $Y$ is irreducible, $s$ must be dominant. Moreover, we may arrange that $X'$ is smooth. Indeed, if it is not smooth, replace $X'$ with a resolution $X''$ and $Y'$ with the fiber product $Y' \times_{X'} X''$. Set $t = f \circ s = r \circ f'$.

Let $Z \subset X$ be a closed subvariety. Proper pushforward gives a commutative diagram of Chow groups
\[
\begin{tikzcd}
CH_0 \big( t^{-1}Z \big) \arrow[r, "s_*"] \arrow[d, "f'_*"'] & CH_0 \big( f^{-1}Z \big) \arrow[d, "f_*"] \\
CH_0 \big( r^{-1}Z \big) \arrow[r, "r_*"'] & CH_0 \big( Z \big).
\end{tikzcd}
\]
Since $X$ has $CH_0$-trivial singularities, $r_*$ is an isomorphism. Moreover, the very general fibers of $f'$ and $f$ are isomorphic. Hence, Theorem \ref{thm flat} below, together with Lemma \ref{lem CH_0 trivial}, implies that the left vertical map $f'_*$ is an isomorphism. It follows from the commutativity of the diagram that $s_*$ is injective. On the other hand, $s_*$ is surjective because $s$ is dominant, by Lemma \ref{lem dominant}. We conclude that $s_*$ is an isomorphism, and consequently $f_*$ is an isomorphism as well.
\end{proof}

It remains to establish the flat case. Here, we do not impose any regularity hypotheses on either $Y$ or $X$.

\begin{theorem}\label{thm flat}
Let $f: Y \rightarrow X$ be a flat projective morphism, with $X$ irreducible. If the very general fiber $F$ has $CH_0(F) \simeq \mathbb{Z}$, then every fiber $F$ has $CH_0(F) \simeq \mathbb{Z}$.
\end{theorem}

\begin{proof}
Fix a closed point $x \in X$. We prove that the fiber $Y_x = f^{-1}(x)$ satisfies $CH_0(Y_x) \simeq \mathbb{Z}$. 

We first reduce to the case where $X$ is a smooth curve. Indeed, if $C \subset X$ is a general curve passing through $x$, we are free to replace $X$ with the normalization $\widetilde{C}$ of $C$ and $Y$ with the fiber product $Y \times_X \widetilde{C}$.

For the remainder of the proof, we write $X = C$, with function field $K = \mathbb{C}(C)$. Let $Y_K$ denote the generic fiber and $Y_{\overline{K}}$ the geometric generic fiber of $f$. By \cite[Lemma $2.1$]{vial}, there is an abstract field isomorphism $\overline{K} \simeq \mathbb{C}$ under which $Y_{\overline{K}}$ is isomorphic as a variety to the very general fiber of $f$. Since a very general fiber satisfies $CH_0(F) \simeq \mathbb{Z}$, it follows that $CH_0(Y_{\overline{K}}) \simeq \mathbb{Z}$.

To prove the isomorphism $CH_0(Y_x) \simeq \mathbb{Z}$, it is enough to show that any two closed points $p$ and $q$ in $Y_x$ have the same class $[p] = [q] \in CH_0(Y_x)$. The strategy is to lift $p$ and $q$, after passing to a suitable finite extension of $K$, to points of the corresponding generic fiber that specialize to $p$ and $q$, respectively. Since $CH_0(Y_{\overline{K}}) \simeq \mathbb{Z}$, these lifted points will become rationally equivalent after a further finite extension. Applying Fulton's specialization map then yields $[p] = [q]$, as required.

We first construct a lift of $p$. Choose a closed point $P$ in the generic fiber $Y_K$ which specializes to $p$. Thus, its closure $W = \overline{ \{ P \} } \subset Y$ is an irreducible curve containing $p$. Let $\widetilde{p}$ be a point in the normalization $\widetilde{W}$ of $W$ lying over $p$, and let $R_p = \mathscr{O}_{\widetilde{W},\widetilde{p}}$ be its local ring. It is a DVR whose fraction field $K_p$ is a finite extension of $K$ and whose residue field is $\mathbb{C}$. The induced morphism $\mathop{\mathrm{Spec}} R_p \rightarrow Y$ sends the closed point to $p$. After base changing $f$ along the induced morphism $\operatorname{Spec}R_p \rightarrow C$, we obtain a flat projective morphism
\begin{align*}
    Y_{R_p}:=Y\times_C\operatorname{Spec}R_p \rightarrow \operatorname{Spec}R_p,
\end{align*}
together with a section $s_p: \operatorname{Spec}R_p\rightarrow Y_{R_p}$ whose special point maps to $p$.

Perform the same construction for $q$, obtaining a morphism $Y_{R_q}: \rightarrow \operatorname{Spec}R_q$, together with a section $s_q: \operatorname{Spec}R_q\rightarrow Y_{R_q}$ whose special point maps to $q$.

Choose a DVR $R$ which dominates both $R_p$ and $R_q$. After a base change to $\operatorname{Spec} R$, we obtain a flat projective morphism
\begin{align*}
    f_R: Y_R \rightarrow \operatorname{Spec} R,
\end{align*}
together with induced sections $s_p$ and $s_q$. The closed fiber is isomorphic to $Y_x$, and the generic fiber is isomorphic to $Y_L$ for a finite extension $L/K$.

After base changing the generic fiber to the algebraic closure $\overline{K}$, the points $P$ and $Q$ determine closed points $P_{\overline{K}}, Q_{\overline{K}} \in Y_{\overline{K}}$. Since $CH_0(Y_{\overline K})\simeq \mathbb Z$, we have $[P_{\overline K}]=[Q_{\overline K}] \in CH_0(Y_{\overline{K}})$. A rational equivalence between these two classes is represented by finitely many curves, along with finitely many rational functions. Hence all of this data is defined over some finite extension $M$ of $L$. After replacing $L$ with $M$ and $R$ with a localization over the closed point of its integral closure in $M$, we may therefore assume that
\begin{align*}
    [P] = [Q] \in CH_0(Y_L).
\end{align*}

We now consider Fulton's specialization homomorphism for the family $f_R$, which is of the form
\begin{align*}
    \sigma_{R}: CH_0(Y_L) \rightarrow CH_0(Y_x);
\end{align*}
see Theorem \ref{thm fulton}. The closure of $P$ in $Y_R$ is the image of the section $s_p$, and its special fiber is the reduced point $p$, with multiplicity one. By the definition of the specialization map, $\sigma_R([P]) = [p]$. Similarly, $\sigma_R([Q]) = [q]$. Applying $\sigma_{R}$ to the equality $[P] = [Q] \in CH_0(Y_L)$ therefore yields $[p] = [q] \in CH_0(Y_x)$.
\end{proof}

The preceding results concern the triviality of the Chow group of zero-cycles over $\mathbb{C}$. For rationality questions, it is useful to consider arbitrary extensions of the base field. Indeed, if $X$ is a smooth complex projective rational variety, then for every field extension $K / \mathbb{C}$, the degree map 
\begin{align*}
    CH_0(X_K) \rightarrow \mathbb{Z}
\end{align*}
is an isomorphism. In this case, we say that $X$ is \emph{universally $CH_0$-trivial}. This obstruction to rationality has found numerous applications in \cite{voisin degen, cp, totaro} by degeneration to singular fibers with very mild singularities. See also \cite{LS} for another application of these ideas involving singular degenerations.

It is natural to ask whether Theorem \ref{intro-thm CH_0} remains valid after replacing ordinary $CH_0$-triviality by this universal notion. The examples below show that this is not the case.

\begin{example}[{\bf Severi-Brauer conic bundle}]\label{exmp enriques conic}
Let $S$ be an Enriques surface, and choose a nonzero class $0 \neq \alpha \in \mathrm{Br}(S)$. Let
\begin{align*}
    f: Y \rightarrow S
\end{align*}
be the associated Severi-Brauer conic. Then $f$ is smooth and projective, and $Y$ is smooth. The general fiber is isomorphic to $\mathbb{P}^1_\mathbb{C}$, which is universally $CH_0$-trivial.

Set $K = \mathbb{C}(S)$. We claim that, nevertheless, the pushforward
\begin{align*}
    f_{K*}: CH_0(Y_K) \rightarrow CH_0(S_K)
\end{align*}
fails to be surjective. For this, consider the natural bilinear pairing
\begin{align*}
    CH_0(S_K) \times \mathrm{Br}(S_K) \rightarrow \mathrm{Br}(K)
\end{align*}
of \cite[Proposition $6.4.2$]{brauer}. Let $[\delta] \in CH_0(S_K)$ be the class of the diagonal, and let $\alpha_K \in \mathrm{Br}(S_K)$ denote the pullback of $\alpha$. Since $S$ is smooth, the natural map $\mathrm{Br}(S) \rightarrow \mathrm{Br}(K)$ is injective. Consequently,
\begin{align*}
    \langle [\delta], \alpha_K \rangle = \delta^* \alpha_K = \alpha_\eta \neq 0 \in \mathrm{Br}(K).
\end{align*}

On the other hand, the defining property of the Severi-Brauer conic gives $f^* \alpha = 0 \in \mathrm{Br}(Y)$, hence $f_K^* \alpha_K = 0$. Therefore, the push-pull formula gives
\begin{align*}
    \langle f_* z, \alpha \rangle = \langle z, f^* \alpha \rangle = 0
\end{align*}
for every $z \in CH_0(Y_K)$. Hence, $[\delta]$ cannot lie in the image of $f_{K*}$.
\end{example}

\begin{example}[{\bf A canonical cone}]\label{exmp klt}
Let $Z = \{ q(x_0, \ldots, x_4) = 0 \} \subset \mathbb{P}^4$ be a very general quartic threefold. By the main result of \cite{cp}, $Z$ is not universally $CH_0$-trivial. That is, for the function field $K = \mathbb{C}(Z)$, the degree map $CH_0(Z_K) \rightarrow \mathbb{Z}$ is not an isomorphism.

Let $X = \{ q(x_0, \ldots ,x_4) = 0 \}\subset \mathbb{P}^5$ be the projective cone over $Z$. Then $X$ has canonical and hence klt singularities. The blow-up $f: Y \rightarrow X$ of $X$ at the vertex is a resolution of singularities. Away from the vertex, $f$ is an isomorphism. Thus the general fiber is a point, so it is certainly universally $CH_0$-trivial. However, after base change to $K = \mathbb{C}(Z)$, the exceptional fiber $E_K$ is isomorphic to $Z_K$, so 
\begin{align*}
    CH_0(E_K) \not\simeq \mathbb{Z}.
\end{align*}
\end{example}

\subsection{Restriction of mixed Hodge modules}\label{subsection restriction}
In this section, we establish general results describing the behavior of mixed Hodge modules and their graded de Rham complexes under a closed embedding
\[
i\colon Z\hookrightarrow X.
\]
Proposition~\ref{prop mhm 2}, concerning the ordinary inverse image $i^*$, is a key ingredient in the proof of Theorem~\ref{intro-thm Omega^p}. On the other hand, Proposition~\ref{prop MHM}, concerning the extraordinary inverse image $i^!$, is needed for the more general form of Theorem~\ref{intro-thm Omega^0} stated as Theorem~\ref{thm Hodge general}. The latter result is not required for the special case presented in the introduction.

We begin with two basic results, analogous to \cite[Lemma~6.1 and Corollary~6.3]{burke lcdef}. In fact, Remark~6.2 therein records the following statement; we include its proof here for the sake of completeness.

\begin{lemma}\label{lem cohom vs grDR}
Suppose $X$ is embedded in a smooth variety $V$. Let $M^\bullet \in D^b (\mathrm{MHM}(X))$, and let $c$ and $m$ be integers. Then 
\begin{align*}
    \mathcal{H}^j \gr_p^F \dr ( M^\bullet )  = 0 \hspace{.4cm} \text{ whenever } p \leq m \text{ and } j+p \leq c
\end{align*}
if and only if 
\begin{align*}
    \mathcal{H}^j \gr_p^F ( M^\bullet )  = 0 \hspace{.4cm} \text{ whenever } p \leq m \text{ and } j+p \leq c.
\end{align*}
\end{lemma}

\begin{proof}
We argue by induction on $m$. We may therefore assume that $\mathcal{H}^j \gr_p^F ( M^\bullet ) = 0$ whenever $p < m$ and $j+p \leq c$. It suffices to prove
\begin{align*}
    \mathcal{H}^j \gr_m^F \dr ( M^\bullet ) \simeq \mathcal{H}^j \gr_m^F ( M^\bullet) \hspace{.4cm} \text{ whenever } j \leq c - m.
\end{align*}

Set $d = \dim V$, and consider the expression 
\begin{align*}
    \gr_m^F \dr ( M^\bullet) = \bigg( \gr_{m-d}^F M^\bullet \rightarrow \Omega_{V}^1 \otimes_\mathscr{O} \gr_{m-d+1}^F M^\bullet \rightarrow \ldots \rightarrow \Omega_{V}^d \otimes_\mathscr{O} \gr_m^F M^\bullet \bigg).
\end{align*}
There is a convergent spectral sequence with $E_1$ page
\begin{align*}
    E_1^{a,b} = \mathcal{H}^b \bigg( \gr_{m+a}^F M^\bullet \otimes \Omega^{d+a}_V \bigg) \Rightarrow \mathcal{H}^{a+b} \gr_m^F \dr ( M^\bullet).
\end{align*}
By assumption, the terms $E_1^{a,b}$ vanish if $a > 0$ or if $a < 0$ and $a + b \leq c - m$.

Hence, if $j \leq c - m$, the spectral sequence shows that $E_1^{0,j} = E_\infty^{0,j}$ and, in fact,
\begin{align*}
    \mathcal{H}^j \gr_m^F \dr ( M^\bullet ) \simeq \mathcal{H}^j \gr_m^F ( M^\bullet ).
\end{align*}
This proves the claim.
\end{proof}

\begin{corollary}\label{cor cohom vs grDR}
\begin{align*}
    \mathcal{H}^j \gr_p^F \dr ( M^\bullet )  = 0 \hspace{.4cm} \text{ whenever } p \leq m \text{ and } j+p \leq c
\end{align*}
if and only if
\begin{align*}
    \gr_p^F \dr ( \mathcal{H}^j  M^\bullet )  = 0 \hspace{.4cm} \text{ whenever } p \leq m \text{ and } j+p \leq c.
\end{align*}
\end{corollary}

\begin{proof}
Apply Lemma \ref{lem cohom vs grDR} to $M^\bullet$ and $\mathcal{H}^j M^\bullet$ separately, and use the fact that $\gr_p^F$ commutes with cohomology.
\end{proof}

We now discuss the inverse image $i^*$, for which the main result is the following. We will only need it when $m \gg 0$ is large.

\begin{proposition}\label{prop mhm 2}
Let $i: Z \hookrightarrow X$ be a closed embedding and let $M^\bullet \in D^b(\mathrm{MHM}(X))$. If $c$ and $m$ are integers and
\begin{align*}
    \mathcal{H}^j \gr_p^F \dr ( M^\bullet ) = 0 \hspace{.5cm} \text{ whenever } p \leq m \text{ and } j+p \leq c,
\end{align*}
then 
\begin{align*}
    \mathcal{H}^j \gr_p^F \dr ( i^* M^\bullet ) = 0 \hspace{.5cm} \text{ whenever } p \leq m \text{ and } j+p \leq c.
\end{align*}
\end{proposition}

\begin{proof}
By Corollary \ref{cor cohom vs grDR}, the hypothesis is equivalent to
\begin{align*}
    \gr_p^F \dr ( \mathcal{H}^j M^\bullet ) = 0 \hspace{.5cm} \text{ whenever } p \leq m \text{ and } j+p \leq c.
\end{align*}

Consider now the spectral sequence
\begin{align*}
    E_2^{a,b} = \gr_p^F \dr ( \mathcal{H}^a i^* \mathcal{H}^b M^\bullet ) \Rightarrow \gr_p^F \dr ( \mathcal{H}^{a+b} i^* M^\bullet ).
\end{align*}
Recall that $i^*$ is right exact for the perverse t-structure. Therefore, applying Lemma \ref{lem i^*} below to each $\mathcal{H}^b M^\bullet$, we obtain vanishing $E_2^{a,b} = 0$ whenever $p \leq m$ and $p \leq c-b-a$. Then the spectral sequence shows that 
\begin{align*}
    \gr_p^F \dr ( \mathcal{H}^{a+b} i^* M^\bullet ) = 0 \hspace{.5cm} \text{ whenever } p \leq m \text{ and } a+b+p \leq c.
\end{align*}
The conclusion follows from another application of Corollary \ref{cor cohom vs grDR}.
\end{proof}

The proof of Proposition~\ref{prop mhm 2} relies on the following result concerning the inverse image of a single mixed Hodge module.

\begin{lemma}\label{lem i^*}
Let $i : Z \rightarrow X$ be a closed embedding and let $M \in \mathrm{MHM}(X)$. If $\gr_p^F \dr (M) = 0$ for all $p \leq m$, then for any $j \leq 0$, we have $\gr_p^F \dr ( \mathcal{H}^{j} i^* M) = 0$ for all $p \leq m-j$.
\end{lemma}

\begin{proof}
This is a local statement, so we may assume $X$ is affine. Moreover, replacing $X$ with a smooth ambient variety if necessary, we may assume it is smooth. Let $( \mathcal{M}, F_\bullet)$ denote the filtered $\mathcal{D}_X$-module underlying $M$. Examining equation \eqref{eqn grDR}, we see that our hypothesis is equivalent to the vanishing $F_m \mathcal{M} = 0$. 

We may further reduce to the case where $Z = \{ g = 0\}$ is defined by a single equation. Indeed, we argue by induction on the number of equations. Suppose that there is a factorization $i = i_1 \circ i_2$, where the lemma is known for both $i_1$ and $i_2$. Consider the spectral sequence
\begin{align*}
    E_2^{a,b} = \mathcal{H}^a i_2^* ( \mathcal{H}^b i_1^* M ) \Rightarrow \mathcal{H}^{a+b} i^* M.
\end{align*}
By the lemma applied to $i_1$, we have $\gr_p^F \dr ( \mathcal{H}^b i_1^* M) = 0$ for $p \leq m - b$. But then the lemma applied to $i_2$ gives $\gr_p^F \dr (E_2^{a,b}) =0$ for $p \leq m -a-b$. Our result then follows from the spectral sequence.

It therefore remains to consider the case in which $Z = \{ g = 0 \}$. Then \cite[Corollary 2.24]{saito 90} shows that $\mathcal{H}^j i^* M = 0$ for $j \neq 0, -1$, and there is an exact sequence
\begin{align*}
    0 \rightarrow i_* \mathcal{H}^{-1} i^* M \rightarrow \psi_{g,1} M \xrightarrow{\mathrm{can}} \phi_{g,1} M \rightarrow i_* \mathcal{H}^0 i^* M \rightarrow 0.
\end{align*}
Consider now the graph embedding $i_g: X \rightarrow X \times \mathbb{C}$, $x \mapsto \big( x, g(x) \big)$ and the pushforward $\widetilde{\mathcal{M}} = (i_g)_+ \mathcal{M}$. Then \cite[Equation $2.2.6$]{saito 90} gives identifications
\begin{align*}
    \psi_{g,1} M = \bigg( \gr_V^1 \widetilde{\mathcal{M}}, F_{\bullet - 1} \gr_V^1 \widetilde{\mathcal{M}} \bigg), \qquad \phi_{g,1} M = \bigg( \gr_V^0 \widetilde{\mathcal{M}}, F_{\bullet} \gr_V^0 \widetilde{\mathcal{M}} \bigg).
\end{align*}
It follows that $F_{m} \phi_{g,1} M = F_{m} \gr_V^0 \widetilde{\mathcal{M}} = 0$, hence the sequence shows that $F_m i_* \mathcal{H}^0 i^* M = 0$. Similarly, $F_{m+1} \psi_{g,1} M = F_{m} \gr_V^1 \widetilde{\mathcal{M}} = 0$, so we conclude $F_{m+1} i_* \mathcal{H}^{-1} i^* M = 0$, as required.
\end{proof}

There is also a converse to Proposition \ref{prop mhm 2} in the more restrictive setting where $m \gg 0$ is large.

\begin{proposition}\label{prop mhm 2 converse}
Let $M^\bullet \in D^b ( \mathrm{MHM}(X))$ and let $c$ be an integer. Suppose that, for every inclusion of a closed point $i_x: \{ x \} \rightarrow X$, one has
\begin{align*}
    \mathcal{H}^j \gr_p^F \dr ( i_x^* M^\bullet ) = 0 \hspace{.4cm} \text{ whenever } j + p \leq c.
\end{align*}
Then
\begin{align*}
    \mathcal{H}^j \gr_p^F \dr ( M^\bullet ) = 0 \hspace{.4cm} \text{ whenever } j + p \leq c.
\end{align*}
\end{proposition}

\begin{proof}
We follow the same argument as in \cite[Lemma $1.8$]{dor}, but using the first distinguished triangle in Proposition \ref{prop dist tri} in place of the second. This reduces us to the case in which $X$ is smooth and each cohomology module $\mathcal{H}^j M^\bullet$ is a variation of mixed Hodge structures.

Let $x \in X$. Since $i_x$ is non-characteristic with respect to each cohomology module $\mathcal{H}^j M^\bullet$, our indexing conventions give
\begin{align*}
    \mathcal{H}^j \gr_p^F \dr ( i_x^* M^\bullet ) \simeq i_x^* \mathcal{H}^{j+n} \gr_{p-n}^F ( M^\bullet),
\end{align*}
where $n= \dim X$.

By assumption, the left hand side vanishes whenever $j + p \leq c$. Hence, 
\begin{align*}
    \mathcal{H}^{j+n} \gr_{p-n}^F ( M^\bullet) = 0 \hspace{.4cm} \text{ whenever } j + p \leq c.
\end{align*}
The result now follows from Lemma~\ref{lem cohom vs grDR}.
\end{proof}

We now turn to the extraordinary inverse image $i^!$. Although somewhat awkward, the following proposition is precisely what is needed to prove Theorem~\ref{thm Hodge general}.

\begin{proposition}\label{prop MHM}
Let $i : Z \rightarrow X$ be a closed embedding, where $n = \dim X$, and let $M^\bullet \in D^b(\mathrm{MHM}(X))$. Suppose that $c$ and $m$ are integers and
\begin{align*}
    \mathcal{H}^j \gr_p^F \dr( M^\bullet) = 0 \hspace{.5cm} \text{ for all } p \leq m \text{ and } j \geq c.
\end{align*}
Then
\begin{align*}
    \mathcal{H}^j \gr_p^F \dr( \mathbb{D} i^! M^\bullet) = 0 \hspace{.5cm} \text{ for all } p \geq -m \text{ and } j \leq -c-n
\end{align*}
\end{proposition}

When $c \ll 0$, Proposition \ref{prop MHM} specializes to \cite[Lemma $1.8$]{dor}. This special case is all that is needed for Theorem \ref{intro-thm Omega^0} stated in the introduction.

\begin{proof}
Fix $p \geq -m$, and consider the spectral sequence
\begin{align*}
    E_2^{a,b} = \mathcal{H}^b \gr_p^F \dr ( i^* \mathcal{H}^a \mathbb{D} M^\bullet ) \Rightarrow \mathcal{H}^{a+b} \gr_p^F \dr ( i^* \mathbb{D} M^\bullet ).
\end{align*}
It is enough to prove that $E_2^{a,b} = 0$ whenever $a + b \leq -c -n$. By Lemma \ref{lem dim vanish} below applied to the mixed Hodge module $\mathcal{H}^a \mathbb{D} M^\bullet$, we have $E_2^{a,b} = 0$ for $b < -n$. It therefore suffices to show that $E_2^{a,b} = 0$ for $a \leq -c$.

Fix $a \leq -c$ as above. By \cite[Corollary $6.3$]{burke lcdef}, the hypothesis on $M^\bullet$ gives
\begin{align*}
    \gr_{q}^F \dr ( \mathcal{H}^{-a} M^\bullet) = 0 \hspace{.5cm} \text{ for all } q \leq m.
\end{align*}
By \cite[Lemma $1.8$]{dor}, this implies
\begin{align*}
    \gr_{q}^F \dr ( i^! \mathcal{H}^{-a} M^\bullet ) = 0 \hspace{.5cm} \text{ for all } q \leq m.
\end{align*}
Set $q = -p$. Taking the Grothendieck dual and applying Proposition \ref{prop D comp} gives
\begin{align*}
    \gr_p^F \dr \big( i^* \mathcal{H}^a \mathbb{D} M^\bullet \big) = 0,
\end{align*}
as desired.
\end{proof}

We required the following lemma for the proof of Proposition \ref{prop MHM} above.

\begin{lemma}\label{lem dim vanish}
Let $i: Z \hookrightarrow X$ be a closed embedding with $n = \dim X$, and let $M \in \mathrm{MHM}(X)$. Then, for every $p \in \mathbb{Z}$,
\begin{align*}
    \mathcal{H}^b \gr_p^F \dr ( i^* M ) = 0 \hspace{.5cm} \text{ for all } b < -n.
\end{align*}
\end{lemma}

\begin{proof}
We first prove that
\begin{align}\label{eqn dim supp}
    \dim \Supp \mathcal{H}^{s} i^*M\leq n+s \hspace{.5cm} \text{ for every } s \in \mathbb{Z}.
\end{align}

Let $P$ be the perverse sheaf underlying $M$. For any integer $s \in \mathbb{Z}$, set $d = \dim \Supp {}^p \! \mathcal{H}^s i^* P$, and let $x$ be a general point in this support. Let $N$ be a normal slice through $x$ of codimension $d$, transverse to a stratification adapted to $P$, $Z$, and $\Supp {}^p \! \mathcal{H}^s i^* P$.

By the non-characteristic restriction theorem, $P_N := P|_N[-d]$ is perverse on $N$, and
\begin{align}\label{eqn neq zero}
    {}^p\!\mathcal{H}^{s} ( i_N^* P_N) \simeq {}^p\!\mathcal{H}^{s}(i^* P)|_{Z \cap N}[-d] \neq 0.
\end{align}

On the other hand, since $P_N$ is perverse and $i_N^*$ is left exact for constructible sheaves, we have $i_N^* P_N \in {}^c\!D^{ \geq -(n-d)}$. But \cite[Proposition $8.1.22$, $8.1.24$]{htt} shows there is an inclusion
\begin{align*}
    {}^c\!D^{ \geq -(n-d)} \subset {}^p\!D^{ \geq -(n-d)}.
\end{align*}
Hence, \eqref{eqn neq zero} implies that $s \geq -(n-d)$. Rearranging, we obtain \eqref{eqn dim supp}.

Consider now the spectral sequence
\begin{align*}
    E_2^{r,s} = \mathcal{H}^r \gr_p^F \dr ( \mathcal{H}^s i^* M ) \Rightarrow \mathcal{H}^{r+s} \gr_p^F \dr ( i^* M).
\end{align*}

By \eqref{eqn dim supp}, we obtain vanishing $E_2^{r,s} = 0$ whenever $r < -n-s$. The spectral sequence shows that $\mathcal{H}^{r+s} \gr_p^F \dr ( i^* M) = 0$ for $r+s < -n$.
\end{proof}

\subsection{Rigidity of holomorphic forms}

The goal of this section is to prove the following rigidity theorem for the low degree cohomology sheaves of the Du Bois complex.

\begin{theorem}\label{thm Hodge singular}
Let $f: Y \rightarrow X$ be a dominant projective morphism, where $Y$ is irreducible and $X$ satisfies condition $E_0$. Suppose that the general fiber $F$ is connected and satisfies $\mathbb{H}^0(F, \underline{\Omega}_F^p) = 0$ for every $p > 0$. Then for any closed subvariety $Z \subset X$ and $p \geq 0$, the natural morphism
\begin{align*}
    \mathcal{H}^0 \underline{\Omega}_Z^p \rightarrow \mathcal{H}^0 Rf_* \underline{\Omega}_{f^{-1}(Z)}^p
\end{align*}
is an isomorphism, and the morphism
\begin{align*}
    \mathcal{H}^1 \underline{\Omega}_Z^p \rightarrow \mathcal{H}^1 Rf_* \underline{\Omega}_{f^{-1}(Z)}^p
\end{align*}
is injective.
\end{theorem}

Passing to hypercohomology gives the following global form of the theorem. 

\begin{corollary}[{\(\!\Rightarrow\) Theorem~\ref{intro-thm Omega^p}}]
Under the hypotheses of Theorem \ref{thm Hodge singular}, for any closed subvariety $Z \subset X$ and $p \geq 0$, the induced morphism
\begin{align*}
    \mathbb{H}^0 \big( Z, \underline{\Omega}_Z^p \big) \rightarrow \mathbb{H}^0 \big( f^{-1}(Z), \underline{\Omega}_{f^{-1}(Z)}^p \big)
\end{align*}
is an isomorphism, and the morphism
\begin{align*}
    \mathbb{H}^1 \big( Z, \underline{\Omega}_Z^p \big) \rightarrow \mathbb{H}^1 \big( f^{-1}(Z), \underline{\Omega}_{f^{-1}(Z)}^p \big)
\end{align*}
is injective. In particular, every fiber $F$ satisfies $\mathbb{H}^0(F, \underline{\Omega}_F^p) = 0$ for $p > 0$.
\end{corollary}

Theorem \ref{intro-thm Omega^p} follows because rational singularities imply condition $D_0$ and hence also condition $E_0$ by Proposition \ref{prop hiatt}.

\begin{remark}[{\bf Differential forms in the h-topology}]
By \cite[Theorem $7.12$]{hj}, there is a natural identification $\mathbb{H}^0(Z, \underline{\Omega}_Z^p) \simeq \Omega_h^p(Z)$, where $\Omega_h^p(X)$ denotes the space of holomorphic $p$-forms on $X$ in the h-topology. 

Moreover, \cite[Lemma $2.9$]{hj} shows that every dominant morphism $f: Y \rightarrow X$ induces an injective pullback map $\Omega_h^p(X) \rightarrow \Omega_h^p(Y)$. This is the Hodge-theoretic counterpart of the surjectivity of the pushforward map $CH_0(Y) \rightarrow CH_0(X)$ established in Lemma \ref{lem dominant}. As in Remark \ref{rmk Y irred}, this perspective also explains why no assumptions on the singularities of $Y$ are needed above.
\end{remark}

We now turn to the proof of Theorem \ref{thm Hodge singular}.

\begin{proof}[Proof of Theorem \ref{thm Hodge singular}]
There is an induced morphism $\mathbb{Q}_X^H \rightarrow f_* \mathbb{Q}_Y^H$ in the derived category $D^b(\mathrm{MHM}(X))$ of mixed Hodge modules. Let $M_f^\bullet$ be its cone, so there is a distinguished triangle
\begin{align}\label{eqn tri M_f}
    \mathbb{Q}_X^H \rightarrow f_* \mathbb{Q}_Y^H \rightarrow M_f^\bullet \xrightarrow[]{+1}
\end{align}
We claim that there is vanishing
\begin{align}\label{eqn M_f vanish}
    \mathcal{H}^j \gr_{-p}^F \dr ( M_f^\bullet) = 0 \hspace{.5cm} \text{ whenever } j-p \leq 0.
\end{align}
Indeed, fix an integer $p$. Applying $\gr_{-p}^F \dr(-)$ to the triangle \eqref{eqn tri M_f} and taking cohomology gives a long exact sequence
\begin{align*}
    0 \rightarrow \mathcal{H}^{p-1} \gr_{-p}^F \dr ( M_f^\bullet ) \rightarrow \mathcal{H}^0 \underline{\Omega}_X^p \xrightarrow{\phi_0} &\mathcal{H}^0 Rf_* \underline{\Omega}_Y^p \rightarrow \mathcal{H}^{p} \gr_{-p}^F \dr ( M_f^\bullet ) \\ &\rightarrow \mathcal{H}^1 \underline{\Omega}_X^p \xrightarrow{\phi_1} \mathcal{H}^1 Rf_* \underline{\Omega}_Y^p \rightarrow \ldots
\end{align*}
The desired vanishing is immediate for $j < p-1$. It therefore remains to prove that $\phi_0$ is an isomorphism and that $\phi_1$ is injective. Since $\mathbb{H}^0(F, \underline{\Omega}_F^q) = 0$ for all $0 < q \leq p$, the morphism $\phi_0$ is generically an isomorphism.

Set $n = \dim X$ and $m = \dim Y$. By the decomposition theorem, $IC_X^H[-n]$ is a direct summand of $f_* IC_Y^H[-m]$. Choose a projection onto this summand, and consider the composition 
\begin{align*}
    \mathbb{Q}_X^H \rightarrow f_* \mathbb{Q}_Y^H \rightarrow f_* IC_Y^H[-m] \rightarrow IC_X^H[-n].
\end{align*}
It is straightforward to see it agrees with the natural map in \eqref{eqn K}. Applying $\gr_{-p}^F \dr (-)$ and taking cohomology gives a commutative diagram
\[
\begin{tikzcd}
\mathcal{H}^0\underline{\Omega}_X^p
    \arrow[r, "\phi_0"]
    \arrow[dr]
&
\mathcal{H}^0 Rf_*\underline{\Omega}_Y^p
    \arrow[d, "\psi_0"]
\\
&
\mathcal{H}^0 I\underline{\Omega}_X^p
\end{tikzcd}
\]

Since $X$ satisfies condition $E_0$, the diagonal map is an isomorphism. Consequently, $\phi_0$ is a split injection. On the other hand, $\mathcal{H}^0 \underline{\Omega}_Y^p$ is torsion-free; see for instance \cite[\S 1]{parkpopa}. Hence, its direct image $f_* \mathcal{H}^0 \underline{\Omega}_Y^p \simeq \mathcal{H}^0 Rf_* \underline{\Omega}_Y^p$ is also torsion-free. Because $\phi_o$ is split, its cokernel is a direct summand of $\mathcal{H}^0 Rf_* \underline{\Omega}_Y^p$ and is therefore torsion-free as well. But we have seen that $\phi_0$ is generically an isomorphism, so this cokernel is torsion. It follows that $\mathrm{coker} \phi_0 = 0$, and hence $\phi_0$ is an isomorphism.

Similarly, we obtain in degree one a commutative diagram
\[
\begin{tikzcd}
\mathcal{H}^1\underline{\Omega}_X^p
    \arrow[r, "\phi_1"]
    \arrow[dr]
&
\mathcal{H}^1 Rf_*\underline{\Omega}_Y^p
    \arrow[d, "\psi_1"]
\\
&
\mathcal{H}^1 I\underline{\Omega}_X^p
\end{tikzcd}
\]
Since $X$ satisfies condition $E_0$, the diagonal arrow is injective. Thus, $\phi_1$ is injective. This implies the claimed vanishing \eqref{eqn M_f vanish}.

Now let $i : Z \rightarrow X$ be a closed embedding. By Proposition \ref{prop mhm 2}, the vanishing \eqref{eqn M_f vanish} propagates to give 
\begin{align*}
    \mathcal{H}^j \gr_{-p}^F \dr ( i^* M_f^\bullet) = 0 \hspace{.5cm} \text{ whenever } j-p \leq 0.
\end{align*}

Apply the exact functor $\gr_{-p}^F \dr (i^*-)[p]$ to the triangle \eqref{eqn tri M_f}. Using Propositions \ref{prop base change} and \ref{prop D comp}, we obtain a distinguished triangle
\begin{align*}
    \underline{\Omega}_Z^p \rightarrow Rf_* \underline{\Omega}_{f^{-1}(Z)}^p \rightarrow \gr_{-p}^F \dr ( i^* M_f^\bullet)[p] \xrightarrow[]{+1}
\end{align*}
Taking its cohomology long exact sequence concludes the proof.
\end{proof}

\begin{remark}[{\bf Condition $E_0$}]\label{rmk E_0 char}
Applying the argument from the proof above to a resolution of singularities $f: Y \rightarrow X$, we find that $X$ satisfies condition $E_0$ if and only if
\begin{align*}
    \mathcal{H}^j \gr_{-p}^F \dr ( M_f^\bullet) = 0 \hspace{.5cm} \text{ whenever } j-p \leq 0.
\end{align*}
Indeed, by \cite[Lemma 3.4]{parkpopa} the morphism
\begin{align*}
    \psi_0: f_* \Omega_Y^p \simeq \mathcal{H}^0 I \underline{\Omega}_X^p
\end{align*}
is always an isomorphism. Moreover, as observed above, $\psi_1$ is a projection. Thus, $\phi_1$ is injective if and only if $\psi_1 \circ \phi_1$ is injective.
\end{remark}

\subsection{Rigidity of the zeroth Du Bois complex}

In this section, we study rigidity phenomena for the zeroth Du Bois complex and prove Theorem \ref{intro-thm Omega^0}. As a warm up to the general statement thereof, we give a more direct proof of the following result.

\begin{proposition}\label{prop ratl sing}
Let $f: Y \rightarrow X$ be a projective morphism, and suppose that $R^p f_* \mathscr{O}_Y = 0$ for some $p > 0$. Then every fiber $F = f^{-1}(x)$ satisfies 
\begin{align*}
    \mathbb{H}^p(F, \underline{\Omega}_F^0) = 0.
\end{align*}
\end{proposition}

The proof is essentially the same as that of \cite[Lemma $2.5$]{ms}, where the vanishing $H^p(F, \mathscr{O}_F) = 0$ is established for a simple normal crossings divisor $F$. In this setting, the two statements coincide, since simple normal crossings divisors have Du Bois singularities, so there is a quasi-isomorphism $\mathscr{O}_F \simeq \underline{\Omega}_F^0$.

\begin{proof}
We work in the analytic topology. Choose a small Stein neighborhood $U$ of $x \in X$, such that the inclusion $j: F \hookrightarrow f^{-1}(U)$ is a homotopy equivalence. Consider the commutative diagram
\[
\begin{tikzcd}[column sep=large, row sep=large]
H^p\bigl(f^{-1}(U),\mathbb{C}\bigr)
  \arrow[r]
  \arrow[d,"j^*"']
&
H^p\bigl(f^{-1}(U),\mathscr O_{f^{-1}(U)}\bigr)
  \arrow[d,"j^*"]
\\
H^p(F,\mathbb{C})
  \arrow[r]
&
H^p(F,\mathscr O_F)
  \arrow[r]
&
\mathbb H^p\bigl(F,\underline{\Omega}_F^0\bigr)
\end{tikzcd}
\]

By GAGA, we have $R^p f^{\mathrm{an}}_* \mathscr{O}^{\mathrm{an}}_Y = 0$. Since $U$ is Stein, we deduce $H^p(f^{-1}(U), \mathscr{O}_{f^{-1}(U)}) = 0$. Therefore, the total composition $H^p(f^{-1}(U), \mathbb{C}) \rightarrow \mathbb{H}^p(F, \underline{\Omega}_F^0)$ is zero.

On the other hand, this total composition is surjective. Indeed, the left vertical arrow $j^*$ is an isomorphism because $j$ is a homotopy equivalence. Moreover, degeneration of the Hodge-de Rham spectral sequence \eqref{Hodge-dR ss} implies that the bottom composition $H^p(F, \mathbb{C}) \rightarrow \mathbb{H}^p(F, \underline{\Omega}_F^0)$ is surjective. We conclude that $\mathbb{H}^p(F, \underline{\Omega}_F^0) = 0$.
\end{proof}

As noted in the introduction, the proposition applies when $f: Y \rightarrow X$ is a resolution of a variety $X$ with rational singularities. Hence, every fiber $F = f^{-1}(x)$ satisfies $\mathbb{H}^p(F, \underline{\Omega}_F^0) = 0$ for $p > 0$.

It is tempting to try applying cohomology and base change directly to the quasi-isomorphism
\begin{align*}
    Rf_* \mathscr{O}_Y \simeq \mathscr{O}_X,
\end{align*}
and instead conclude that $H^p(F, \mathscr{O}_F) = 0$. This conclusion is false in general. The question of whether such a vanishing holds was raised in the discussion following \cite[Lemma $2.5$]{ms}. The next two examples show that it may fail both for the scheme-theoretic fiber $F = f^{-1}(x)$ and for its reduction $F_{\mathrm{red}}$.

\begin{example}[{\bf Rational singularity with $H^1(F, \mathscr{O}_F) \neq 0$}]\label{exmp ratl sing}
We construct a log resolution $f: Y \rightarrow X$ of a fivefold $X$ with rational singularities such that some scheme-theoretic fiber $F$ satisfies $H^1(F, \mathscr{O}_F) \neq 0$.

Let $S$ be the blow-up of $\mathbb{P}^2$ at eight points in general position, so that $S$ is a degree one del Pezzo surface. The anticanonical divisor $-K_S$ is ample, and the complete linear system $|-K_S|$ is a pencil with a unique reduced basepoint $p$; see for instance \cite[Example $2.12$]{hv}. In addition, $H^i(S, (-K_S)^m) = 0$ whenever $i > 0$ and $m \geq 0$ by Kodaira vanishing, and $-2K_S$ is basepoint-free.

Set $E = S \times S$ and $L = -K_S \boxtimes -K_S = -K_E$, an ample line bundle. Consider the affine cone 
\begin{align*}
    X = \mathrm{Cone}(E, L) = \mathop{\mathrm{Spec}} \big( \bigoplus_{m \geq 0} H^0(E, L^m) \big).
\end{align*} 
It has rational singularities because of the cohomological vanishing $H^i(E, L^m) = 0$ whenever $i > 0$ and $m \geq 0$; see \cite[Proposition B.$3$]{ks}. Blowing up the vertex $v$ of the cone yields a log resolution
\begin{align*}
    f: Y = \mathrm{Tot}(L^{-1}) = \mathop{\mathrm{Spec}}_E \bigg( \bigoplus_{m \geq 0} L^m \bigg) \rightarrow X.
\end{align*}
Thus, $Y$ is smooth, $f$ is an isomorphism away from the vertex $v \in X$, and the scheme-theoretic fiber $F = f^{-1}(v)$ has smooth reduction $F_{\mathrm{red}} \simeq E$. However, $F$ is not reduced. It is this non-reduced structure which we exploit.

Viewing $Y$ as the total space of $L^{-1}$, consider the affine projection $\pi: Y \rightarrow E$. The fiber $F$ is defined by the graded subsheaf of ideals $\bigoplus_{m \geq 1} I_m$ of $\bigoplus_{m \geq 0} L^m$, where
\begin{align*}
    I_m= \mathrm{im} \bigg( \bigoplus_{i = 1}^m H^0(E, L^i) \otimes_\mathbb{C} L^{m-i} \rightarrow L^{m} \bigg).
\end{align*}
Since $L^2$ is basepoint-free, we have $I_m = L^m$ for all $m \geq 2$. Moreover, $I_1 = L \otimes \mathscr{I}_Z$, where $Z = (\{ p \} \times S) \cup (S \times \{p\} )$ is the base scheme of the linear system $|L|$. 

Therefore, 
\begin{align*}
   \pi_* \mathscr{O}_F = \bigoplus_{m \geq 0} L^m/I_m = \mathscr{O}_E \oplus (L \otimes \mathscr{O}_Z).
\end{align*}
Hence, since $\pi$ is affine,
\begin{align*}
    H^1(F, \mathscr{O}_F) = H^1(E, \mathscr{O}_E) \oplus H^1(Z, L|_Z).
\end{align*}
It is easy to see that $H^1(E, \mathscr{O}_E) = 0$, but we will prove that $H^1(Z, L|_Z) = \mathbb{C}$. For this, consider the Mayer-Vietoris sequence
\begin{align*}
    0 \rightarrow L|_Z \rightarrow L|_{\{p \} \times S} \oplus L|_{ S \times \{p \} } \rightarrow L_{ \{p\} \times \{p\}} \rightarrow 0,
\end{align*}
with induced long exact sequence of cohomology groups
\begin{align*}
    \ldots \rightarrow H^0(S, -K_S) \oplus H^0(S, -K_S) \xrightarrow{\delta} \mathbb{C} \rightarrow H^1(Z, L|_Z) \rightarrow H^1(S, -K_S) \oplus H^1(S, -K_S) \rightarrow \ldots 
\end{align*}
We have $H^1(S, -K_S) = 0$, and the connecting map $\delta$ is zero because $p$ is a basepoint of $|-K_S|$. The conclusion follows.
\end{example}

\begin{example}[{\bf Smooth base with $H^1(F_{\mathrm{red}}, \mathscr{O}_{F_{\mathrm{red}}}) \neq 0$}]\label{exmp reduced}
Let $X = \mathbb{A}^3$. We construct a composition $f$ of a sequence of smooth blow-ups
\begin{align*}
    Y = Y_5 \xrightarrow[]{f_5} Y_4 \xrightarrow[]{f_4} Y_3 \xrightarrow[]{f_3} Y_2 \xrightarrow[]{f_2} Y_1 \xrightarrow[]{f_1} X = \mathbb{A}^3.
\end{align*}
We will denote by $D_i = \mathrm{Ex}(f_i)$ the exceptional divisor of $f_i$ and, by abuse of notation, its strict transform on any higher blow-up as well.

The first map $f_1$ is the blow-up of the origin $0 \in \mathbb{A}^3$. The second map $f_2$ is the blow-up of a line in the exceptional plane $D_1$. Next, on $Y_2$, we blow up a conic in $D_1$ which is tangent to the line $D_1 \cap D_2$. On $Y_3$, blow up the unique intersection point in $D_1 \cap D_2 \cap D_3$.

Let us consider $Y_4$. A straightforward computation shows the union $E = D_1 \cup D_2 \cup D_3 \cup D_4$ is given in local coordinates by the union of four planes $xyz(x+y) = 0$ meeting at a point $p$. We finally take $Y_5$ to be the blow-up of $p$ in $Y_4$. We claim that the scheme-theoretic fiber $F = f^{-1}(0)$ has nonvanishing cohomology group $H^1(F_{\mathrm{red}}, \mathscr{O}_{F_\mathrm{red}}) \simeq \mathbb{C}$.

For this, consider the distinguished triangle obtained by pushing forward along $f_5$ the short exact sequence defining $\mathscr{O}_{F_{\mathrm{red}}}$:

\begin{align}\label{eqn tri L}
    Rf_{5*} \mathscr{O}_Y(-{F_{\mathrm{red}}}) \rightarrow Rf_{5*} \mathscr{O}_{Y} \rightarrow Rf_{5*} \mathscr{O}_{F_{\mathrm{red}}} \xrightarrow[]{+1}.
\end{align}
Since $f_5$ is the blow-up of a smooth center, $Rf_{5*} \mathscr{O}_Y \simeq \mathscr{O}_{Y_4}$. We now compute $Rf_{5*} \mathscr{O}_Y(-{F_{\mathrm{red}}})$. Letting $D = D_5 \subset Y$, the key is the pullback identity
\begin{align*}
    f_5^* E = \tilde{E} + 4 D = F_{\mathrm{red}} + 3D.
\end{align*}
It implies $R^bf_{5*} \mathscr{O}_Y(-{F_{\mathrm{red}}}) \simeq \mathscr{O}_{Y_4}(-E) \otimes R^b f_{5*} \mathscr{O}_Y(3D)$ for all $b$. 

Since $D \simeq \mathbb{P}^2$, there are short exact sequences 
\begin{align*}
    0 \rightarrow \mathscr{O}_Y((k-1)D) \rightarrow \mathscr{O}_Y(kD) \rightarrow \mathscr{O}_{\mathbb{P}^2}(-k) \rightarrow 0
\end{align*}
for $k = 1,2,3$. Pushing them forward along $f_{5}$, we see that $f_{5*} \mathscr{O}_Y(3D) = \mathscr{O}_{Y_4}$, $R^1f_{5*} \mathscr{O}_Y(3D) = 0$, and 
\begin{align*}
    R^2 f_{5*} \mathscr{O}_Y(3D) = H^2(\mathbb{P}^2, \mathscr{O}_{\mathbb{P}^2}(-3))_p = \mathbb{C}_p,
\end{align*}
supported at the blow-up center $p \in Y_4$. Putting this all together in the distinguished triangle \eqref{eqn tri L}, we see $f_{5*} \mathscr{O}_{F_{\mathrm{red}}} \simeq \mathscr{O}_E$, $R^1 f_{5*} \mathscr{O}_{F_{\mathrm{red}}} \simeq \mathbb{C}_p$, and $R^2 f_{5*} \mathscr{O}_{F_{\mathrm{red}}} = 0$.

Consider now the Leray spectral sequence
\begin{align*}
    E_2^{a,b} = H^a(E, R^b f_{5*} \mathscr{O}_{F_{\mathrm{red}}} ) \Rightarrow H^{a+b}(F_{\mathrm{red}}, \mathscr{O}_{F_{\mathrm{red}}}).
\end{align*}
For $b = 0$, we have $E_2^{a,0} = H^a(E, \mathscr{O}_E)$. A simple check using the above coordinate chart for $E$ shows these cohomology groups vanish for $a > 0$. Therefore, $E_2^{0,0} = E_2^{0,1} = \mathbb{C}$, and all other terms are zero. It follows that the spectral sequence degenerates and $H^{1}(F_{\mathrm{red}}, \mathscr{O}_{F_{\mathrm{red}}}) = \mathbb{C}$.
\end{example}

Although the preceding examples rule out a general vanishing theorem for the structure sheaf of a fiber, one still obtains vanishing in the top cohomological degrees. This is consistent with what we saw in Examples \ref{exmp ratl sing} and \ref{exmp reduced}.

\begin{proposition}\label{prop pos ratl}
Let $f: Y \rightarrow X$ be a projective morphism, and let $F = f^{-1}(x)$ be a scheme-theoretic fiber of dimension $d$.
\begin{enumerate}
    \item If $R^p f_* \mathscr{O}_Y = 0$ for $p \geq d$, then every closed subscheme $E \subset F$ (including $E= F$ and $E = F_{\mathrm{red}}$), satisfies $H^d(E, \mathscr{O}_E) = 0$. 
    \item If $R^p f_* \mathscr{O}_Y = 0$ for $p \geq d-1$, then $H^{d-1}(F, \mathscr{O}_F) = 0$.
\end{enumerate}
\end{proposition}

Part $(1)$ was proved by a different argument in \cite[Example $2.6$]{ms}.

\begin{proof}
$(1)$ After shrinking $X$ to an affine neighborhood of $x$, we may assume that $X$ is affine. We first treat the case $E = F$. In this case, there is a hypercohomology spectral sequence of Tor sheaves
\begin{align*}
    E_2^{a,b} = H^a \bigg( Y, \mathcal{T}or_{-b}^{f^{-1}\mathscr{O}_X}\big(\mathscr{O}_Y, f^{-1}k(x) \big) \bigg) \Rightarrow \mathbb{H}^{a+b} \bigg( Y, \mathscr{O}_Y \otimes_{f^{-1}\mathscr{O}_X}^\mathbf{L} f^{-1}k(x) \bigg),
\end{align*}
where $k(x) = \mathscr{O}_X/ \mathrm{m}_x$. For $b = 0$, we have $E_2^{a,0}  = H^a(F, \mathscr{O}_F)$. Moreover, $E_2^{a,b} = 0$ for $b > 0$ since negative Tor sheaves vanish. In fact, all of the Tor sheaves appearing in the spectral sequence are supported on $F$, so Grothendieck vanishing also gives $E_2^{a,b} = 0$ for $a > d = \dim F$. It follows from the spectral sequence that
\begin{align}\label{eqn 1 degen}
    H^d(F, \mathscr{O}_F) = E_2^{d,0} = E_\infty^{d,0} = \mathbb{H}^d \big( Y,\mathscr{O}_Y \otimes_{f^{-1}\mathscr{O}_X}^\mathbf{L} f^{-1}k(x) \big).
\end{align}

Using the derived projection formula, the hypercohomology groups of the ``derived fiber'' $\mathscr{O}_Y \otimes_{f^{-1}\mathscr{O}_X}^\mathbf{L} f^{-1}k(x)$ are given by
\begin{align*}
    R \Gamma \big( Y, \mathscr{O}_Y \otimes_{f^{-1}\mathscr{O}_X}^\mathbf{L} f^{-1}k(x) \big) \simeq Rf_* \big( Lf^*( k(x)) \big) \simeq Rf_* \mathscr{O}_Y \otimes^\mathbf{L}_{\mathscr{O}_X} k(x).
\end{align*}
Our hypothesis on $Rf_* \mathscr{O}_Y$, together with the right exactness of tensor product implies that $\mathbb{H}^d(Y, \mathscr{O}_Y \otimes_{f^{-1}\mathscr{O}_X}^\mathbf{L} f^{-1}k(x)) = 0$. Combining this with \eqref{eqn 1 degen}, we deduce $H^d(F, \mathscr{O}_F) = 0$.

In the general case, there is a surjection $\mathscr{O}_F \rightarrow \mathscr{O}_E$. If $\mathcal{K}$ denotes the kernel, there is an induced long exact sequence
\begin{align*}
    \ldots \rightarrow H^d(F, \mathcal{K}) \rightarrow H^d(F, \mathscr{O}_F) \rightarrow H^d(E, \mathscr{O}_E) \rightarrow H^{d+1}(F, \mathcal{K}) \rightarrow \ldots 
\end{align*}
We have seen $H^d(F, \mathscr{O}_F) = 0$, and $H^{d+1}(F, \mathcal{K}) = 0$ since $\mathcal{K}$ is supported on $F$ of dimension $d$. We conclude that $H^d(E, \mathscr{O}_E)$ = 0.

$(2)$ Returning to the Tor spectral sequence, we find that 
\begin{align*}
    H^{d-1}(F, \mathscr{O}_F) = E_2^{d-1, 0} = E_\infty^{d-1,0}, 
\end{align*}
and this is a quotient of $\mathbb{H}^{d-1}(Y, \mathscr{O}_Y \otimes_{f^{-1}\mathscr{O}_X}^\mathbf{L} f^{-1}k(x))$. The latter vanishes, by the same argument as in $(1)$.
\end{proof}

We now begin the proof of the main result of this section. It implies Theorem \ref{intro-thm Omega^0} by taking $k \gg 0$ sufficiently large.

\begin{theorem}[{\(\!\Rightarrow\) Theorem~\ref{intro-thm Omega^0}}]\label{thm Hodge general}
Let $f: Y \rightarrow X$ be a dominant projective morphism between connected varieties satisfying condition $D_0$ and let $k \geq 0$. Suppose that the general fiber $F$ is connected with $\mathbb{H}^p(F, \underline{\Omega}_F^0) = 0$ for every $0 < p \leq k$. Then for any closed subvariety $Z \subset X$, the induced map
\begin{align*}
     \mathcal{H}^p \underline{\Omega}_Z^0 \rightarrow R^p f_* \underline{\Omega}_{f^{-1}(Z)}^0
\end{align*}
is an isomorphism for $p \leq k$ and injective for $p = k+1$. 

In particular, the map
\begin{align*}
    \mathbb{H}^p(Z, \underline{\Omega}_Z^0) \rightarrow \mathbb{H}^p(f^{-1}(Z), \underline{\Omega}_{f^{-1}(Z)}^0)
\end{align*}
is an isomorphism for $p \leq k$ and injective for $p = k+1$.
\end{theorem}

The key input is the following consequence of Kollár’s results on higher direct images of dualizing sheaves \cite{dualizing}. Saito proved analogous statements for the lowest nonzero graded piece of the Hodge filtration on an arbitrary mixed Hodge module in \cite{saitokollar}, though we will not need it here.

\begin{lemma}\label{lem alpha}
Consider a commutative diagram
\[
\begin{tikzcd}
Y' \arrow[r, "s"] \arrow[d, "f'"'] & Y \arrow[d, "f"] \\
X' \arrow[r, "r"'] & X,
\end{tikzcd}
\]
where $f$ is a dominant projective morphism of connected varieties and $r$ and $s$ are resolutions of singularities. Set $m = \dim Y$, $n = \dim X$, and $t = f \circ s = r \circ f'$. Suppose that the general fiber $F'$ of $f'$ is connected and satisfies $H^p(F', \mathscr{O}_{F'}) = 0$ for every $0 < p \leq k$. Then the natural morphism
\begin{align*}
    \alpha: Rt_* \omega_{Y'}[m] \rightarrow Rr_* \omega_{X'}[n]
\end{align*}
induces an isomorphism of cohomology sheaves in degree $i \geq -n-k$ and a surjection in degree $i = -n-k-1$.
\end{lemma}

\begin{proof}
By the Grauert-Riemenschneider vanishing theorem, the derived pushforward $Rr_* \omega_{X'}$ is quasi-isomorphic to the sheaf $r_* \omega_{X'}$, concentrated in degree zero. 

Moreover, we have vanishing
\begin{align*}
    R^a r_* R^b f'_* \omega_{Y'} = 0
\end{align*}
for all $a > 0$ by \cite[Theorem $3.8$]{dualizing}. Thus, the spectral sequence 
\begin{align*}
    E_2^{a,b} = R^a r_* R^b f'_* \omega_{Y'} \Rightarrow R^{a+b} t_* \omega_{Y'}
\end{align*}
degenerates, giving $R^i t_* \omega_{Y'} \simeq r_* R^i f'_* \omega_{Y'}$ for all $i$.

Therefore, in cohomological degree $-n$, the morphism $\mathcal{H}^{-n}(\alpha)$ is of the form
\begin{align*}
    \mathcal{H}^{-n}(\alpha): r_* R^{m-n} f'_* \omega_{Y'} \rightarrow r_* \omega_{X'}.
\end{align*}
By \cite[Proposition $7.6$]{dualizing}, the natural morphism $R^{m-n} f'_* \omega_{Y'} \rightarrow \omega_{X'}$ is an isomorphism. Consequently, $\mathcal{H}^{-n}(\alpha)$ is an isomorphism as well.

It remains to prove the vanishing of $R^i t_* \omega_{Y'} \simeq r_* R^i f'_* \omega_{Y'}$ for all $m-n-k \leq i < m-n$. In this range, the generic fiber of $R^i f'_* \omega_{Y'}$ vanishes by hypothesis. On the other hand, $R^i f'_* \omega_{Y'}$ is torsion-free by \cite[Theorem $2.1$]{dualizing}, so it must vanish. This concludes the proof.
\end{proof}

We are now prepared to prove Theorem \ref{thm Hodge general}.

\begin{proof}[Proof of Theorem \ref{thm Hodge general}]
Set $m = \dim Y$ and $n = \dim X$. Let $M_f \in D^b( \mathrm{MHM}(X))$ be the cone of the induced morphism $\mathbb{Q}_X^H \rightarrow f_* \mathbb{Q}_X^H$ in $D^b(\mathrm{MHM}(X))$. Dualizing gives a distinguished triangle of the form
\begin{align}\label{eqn M}
    \mathbb{D} M_f \rightarrow \mathbb{D} f_* \mathbb{Q}_Y^H \rightarrow \mathbb{D} \mathbb{Q}_X^H \xrightarrow[]{+1}
\end{align}

We claim that
\begin{align}\label{eqn DM_f vanish}
    \mathcal{H}^j \gr_p^F \dr ( \mathbb{D} M_f) = 0 \hspace{.5cm} \text{ for all } p \leq 0 \text{ and } j \geq -n-k.
\end{align}
This is clear for $p < 0$. 

To treat the case $p = 0$, choose compatible resolutions $r: X' \rightarrow X$ and $s: Y' \rightarrow Y$. This puts us in the situation of Lemma \ref{lem alpha}. Indeed, to verify the vanishing hypothesis for the general fiber $F'$ of $f'$, note that it is a resolution of the general fiber $F$ of $f$. Since $Y$ satisfies condition $D_0$, so too does $F$ by \cite[Remark $4.8(2)$]{dor}. Thus, $\mathbb{H}^p(F, \underline{\Omega}_F^0) \simeq H^p(F' , \mathscr{O}_{F'})$ for all $p \geq 0$.

Since $X$ satisfies condition $D_0$, Lemma \ref{lem D0} and Proposition \ref{prop D comp} give
\begin{align*}
    \gr_0^F \dr \big( \mathbb{D} \mathbb{Q}_X^H \big) \simeq \big(\mathbb{D}_X \underline{\Omega}_X^0 \big)[n] \simeq Rr_* \omega_{X'}[n].
\end{align*}
Similarly, if $t = f \circ s$, we obtain 
\begin{align*}
    \gr_0^F \dr ( \mathbb{D} f_* \mathbb{Q}_Y^H) \simeq Rt_* \omega_{Y'}[m].
\end{align*} 
Now apply $\gr_0^F \dr (-)$ to \eqref{eqn M}. This gives a distinguished triangle
\begin{align*}
    \gr_0^F \dr ( \mathbb{D} M_f ) \rightarrow Rt_* \omega_{Y'}[m] \xrightarrow{\alpha} Rr_* \omega_{X'}[n] \xrightarrow[]{+1}
\end{align*}
By Lemma \ref{lem alpha}, we obtain the desired vanishing \eqref{eqn DM_f vanish}.

Let $i: Z \hookrightarrow X$ be a closed embedding. Applying Proposition \ref{prop MHM}, we obtain vanishing
\begin{align*}
    \mathcal{H}^j \gr_p^F \dr ( i^* M_f)  = 0 \hspace{.5cm} \text{ for all } p \geq 0 \text{ and } j \leq k.
\end{align*}
Finally, we apply $\gr_0^F \dr (i^* \mathbb{D} -)$ to \eqref{eqn M}. This yields a distinguished triangle
\begin{align*}
    \underline{\Omega}_Z^0 \rightarrow Rf_* \underline{\Omega}_{f^{-1}(Z)}^0 \rightarrow \gr_0^F \dr ( i^* M_f) \xrightarrow[]{+1}
\end{align*}
Passing to the associated long exact sequence of cohomology sheaves, the preceding vanishing completes the proof.
\end{proof}

The following example shows that, unlike in Theorem \ref{intro-thm Omega^p}, it is essential that the total space $Y$ satisfy condition $D_0$ in Theorem \ref{intro-thm Omega^0}.

\begin{example}[{\bf Degeneration to a nodal curve}]\label{exmp nodal}
Consider the projection $\mathrm{pr}_1: \mathbb{P}^1 \times \mathbb{P}^1 \rightarrow \mathbb{P}^1 = X$. Fix a point $x_0 \in X$, and identify two distinct points in the fiber $\mathrm{pr}_1^{-1}(x_0)$. Let $Y$ be the resulting variety, and let $f: Y \rightarrow X$ be the induced morphism.

For every $x \neq x_0$, we have $f^{-1}(x) \simeq \mathbb{P}^1$, and so $H^p(f^{-1}(x), \mathscr{O}_{f^{-1}(x)}) = 0$ for $p > 0$. However, the special fiber $F = f^{-1}(x_0)$ is an irreducible nodal rational curve, where
\begin{align*}
    \mathbb{H}^1(F, \underline{\Omega}_F^0) \simeq \mathbb{C}
\end{align*}
is non-vanishing. On the other hand, $\mathbb{H}^0(F, \underline{\Omega}_F^1) = 0$, as predicted by Theorem \ref{intro-thm Omega^p}.
\end{example}

\section{Chow-Hodge parallels}

\subsection{Decomposition of the diagonal}

Recall the generalized Bloch conjecture, Conjecture \ref{intro-conj BB}, predicts that for every smooth, connected, complex projective variety $X$, one has 
\[
CH_0(X)\simeq \mathbb{Z}
\qquad\Longleftrightarrow\qquad
H^0(X,\Omega_X^p)=0
\quad\text{for every }p>0.
\]
The implication from left to right was proved by Roitman \cite{roitman}, and Bloch and Srinivas gave an alternative proof using the decomposition of the diagonal \cite{bloch srinivas}.

In this section, we extend the decomposition of the diagonal argument to singular varieties, replacing $H^0(X, \Omega_X^p)$ with the Hodge-Du Bois piece $\mathbb{H}^0\!\left(X,\underline{\Omega}_X^p\right)$. This is the natural invariant appearing in the rigidity theorems of this paper.

\begin{theorem}[{\(=\) Theorem~\ref{intro-thm BB}}]\label{thm BB}
Let $X$ be a connected projective variety of pure dimension. If the degree map $CH_0(X) \rightarrow \mathbb{Z}$ is an isomorphism, then 
\[
\mathbb{H}^0\!\left(X,\underline{\Omega}_X^p\right)=0
\qquad\text{for every }p>0.
\]
\end{theorem}

We expect the converse to hold, although it remains open even for smooth varieties. Corollary~\ref{cor CH_0 = Z} below provides motivic evidence in support of this expectation.

The proof of Theorem \ref{thm BB} applies a decomposition of the diagonal result for singular varieties from \cite{laterveer} to decompose the cap product with the fundamental class. The key Hodge-theoretic input needed to complete the argument is as follows.

\begin{lemma}\label{lem H^p,0}
Let $X$ be a projective variety of pure dimension $n$. Then cap product with the fundamental class
\begin{align*}
    H^p(X, \mathbb{C}) \xrightarrow[]{\cap [X]} H_{2n-p}(X, \mathbb{C}),
\end{align*}
is injective on the Hodge-Du Bois piece $H^{p,0}(X)
    :=\mathbb{H}^0(X,\underline{\Omega}_X^p)
    \simeq \gr_{-p}^F H^p(X,\mathbb{C})$.
\end{lemma}

\begin{proof}
Given an integer $k$, let
\begin{align*}
    \gamma_k: H^k(X, \mathbb{C}) \rightarrow IH^k(X, \mathbb{C})
\end{align*}
denote the natural comparison morphism, and let $\gamma_{p,q}$ denote its restriction to $H^{p,q}(X)$. Intersection cohomology is pure and it satisfies Poincaré duality. Moreover, cap product map with the fundamental class $\cap [X]$ is the composition
\begin{align*}
    H^p(X, \mathbb{C}) \xrightarrow[]{\gamma_p} IH^p(X, \mathbb{C}) \xrightarrow{\mathrm{PD}_{\mathrm{IH}}} IH^{2n-p}(X, \mathbb{C})^\vee \xrightarrow{\gamma_{2n-p}^\vee} H^{2n-p}(X, \mathbb{C})^\vee \simeq H_{2n-p}(X, \mathbb{C}).
\end{align*}

We restrict this factorization to the Hodge-Du Bois piece $H^{p,0}(X)$. Since the intersection pairing pairs $IH^{p,0}(X)$ perfectly with $IH^{n-p,n}(X)$, we obtain the factorization
\begin{align*}
    H^{p,0}(X) \xrightarrow{\gamma_{p,0}} IH^{p,0}(X) \xrightarrow{\simeq} IH^{n-p, n}(X)^\vee \xrightarrow{\gamma_{n-p,n}^\vee} H^{n-p,n}(X)^\vee.
\end{align*}
It therefore suffices to show that both $\gamma_{p,0}$ and $\gamma_{n-p,n}^\vee$ are injective. For the former, Weber's theorem identifies
\begin{align*}
    \ker \gamma_{p,0} \simeq \gr_{-p}^F W_{p-1} H^p(X, \mathbb{C}),
\end{align*}
and this clearly vanishes; see \cite[Remark $6.5$]{parkpopa}.

For the latter, \cite[Remark $7.2$]{parkpopa} applied with $m = -1$, gives an isomorphism 
\begin{align*}
    \gamma_{n,n-p}: H^{n,n-p}(X) \xrightarrow{\simeq} IH^{n,n-p}(X)
\end{align*}
Moreover, the Hodge-Du Bois piece $H^{n,n-p}(X)$ lies entirely in weight $2n-p$. Taking complex conjugates on this top-weight piece therefore shows that 
\begin{align*}
    \gamma_{n-p,n}: H^{n-p,n}(X) \rightarrow IH^{n-p,n}(X)
\end{align*}
is surjective. We conclude that its dual $\gamma_{n-p,n}^\vee$ is injective, as required.
\end{proof}

We now present the proof of Theorem \ref{thm BB}.

\begin{proof}[Proof of Theorem \ref{thm BB}]
Set $n = \dim X$, and fix a closed point $x \in X$. 

Since $CH_0(X) \simeq \mathbb{Z}$, the decomposition of the diagonal result in \cite[Lemma $1$]{laterveer} gives, after clearing denominators, an identity of the form
\begin{align}\label{eqn decomp}
    N [\Delta_X] = Z_1 + Z_2 \in CH_n(X).
\end{align}
Here, $N$ is a nonzero integer, and $Z_1$ and $Z_2$ are represented by cycles with
\begin{align*}
    \mathrm{Supp} (Z_1) \subset \{ x\} \times X, \hspace{.5cm} \mathrm{Supp} (Z_2) \subset X \times D,
\end{align*}
for some proper subvariety $D \subsetneq X$ of dimension $< n$.

By \cite[Definition $2$]{laterveer}, every correspondence $C \in CH_n(X \times X)$ induces a map
\begin{align*}
    C_*: H^p(X, \mathbb{C}) \rightarrow H_{2n-p}(X, \mathbb{C}), \hspace{.5cm} C_*(\alpha) := (\mathrm{pr}_2)_* \big( \mathrm{pr_1}^*(\alpha) \cap C \big).
\end{align*}
Moreover, \cite[Lemma $4$]{laterveer} shows that the class of the diagonal $\Delta_X$ acts by cap product with the fundamental class.

Fix $p > 0$ and suppose $\alpha \in \mathbb{H}^0(X, \underline{\Omega}_X^p) \simeq \gr_{-p}^F H^p(X, \mathbb{C})$. Applying the correspondence action to \eqref{eqn decomp} gives
\begin{align}\label{eqn sum}
   N \big( \alpha \cap [X] \big) = (Z_1)_*(\alpha) + (Z_2)_*(\alpha). 
\end{align}

We will show that both terms on the right vanish. Indeed, since $Z_1$ is supported on $\{ x \} \times X$, \cite[Lemma $5$]{laterveer} gives a factorization of $(Z_1)_*$ as 
\begin{align*}
    H^p(X, \mathbb{C}) \rightarrow H^p( \{ x \}, \mathbb{C}) \rightarrow H_{2n-p}(X, \mathbb{C}) (-n).
\end{align*}
Because $p > 0$, the middle group vanishes, and hence $(Z_1)_*(\alpha) = 0$.

Similarly, since $Z_2$ is supported on $X \times D$, \cite[Lemma $5$]{laterveer} and compatibility with the Hodge filtration gives a factorization
\begin{align*}
    \gr_{-p}^F H^p(X, \mathbb{C}) \rightarrow \gr_{n-p}^F H_{2n-p}( \widetilde{D}, \mathbb{C}) \rightarrow \gr_{n-p}^F H_{2n-p} (X, \mathbb{C}),
\end{align*}
for a resolution of singularities $\widetilde{D} \rightarrow D$. Since $\dim D < n$, the middle group vanishes. Therefore, $(Z_2)_*(\alpha) = 0$.

Equation \eqref{eqn sum} now gives $\alpha \cap [X] = 0$. But the injectivity result from Lemma \ref{lem H^p,0} therefore implies $\alpha = 0$. Since $\alpha$ was arbitrary, we conclude that $\mathbb{H}^0(X, \underline{\Omega}_X^p) = 0$.
\end{proof}

\subsection{Bloch-Beilinson conjectures for complexes}\label{section BB}
In this section, we formulate motivic heuristics for the relationship between zero-cycles and holomorphic forms on singular varieties; see Corollary \ref{cor CH_0 = Z}. These heuristics will also guide our study of Chow-theoretic singularity conditions in \S \ref{subsection Chow Dm}

To clarify the precise relationship between zero-cycles and holomorphic forms suggested by Conjecture \ref{intro-conj BB}, Bloch and Beilinson predicted that the Chow groups of a smooth projective variety carry a filtration governed by Hodge theory.

\begin{conjecture}\cite[\S $11.2.2$]{voisin ii}\label{conj bb}
Let $X$ be a smooth complex projective variety. Then there exists a decreasing filtration $F^\bullet CH^j(X)_\mathbb{Q}$ on the Chow groups of $X$ with rational coefficients, such that:
\begin{enumerate}
    \item $F^0 CH^j(X)_\mathbb{Q} = CH^j(X)_\mathbb{Q}$, $F^1 CH^j(X)_\mathbb{Q} = CH^j(X)_{\mathrm{hom}, \mathbb{Q}}$, and $F^2 CH^j(X)_\mathbb{Q} = \ker AJ_X$,
    \item $F^{j+1} CH^j(X)_\mathbb{Q} = 0$,
    \item The filtration $F^\bullet CH^j(X)_\mathbb{Q}$ is preserved by correspondences: the map $\Gamma_*$ induced by a cycle $\Gamma \in CH^\ell(X \times Y)$ satisfies
    \begin{align*}
        \Gamma_*(F^\nu CH^j(X)_\mathbb{Q}) \subset F^\nu CH^{j+\ell-\dim X}(Y)_\mathbb{Q},
    \end{align*}
    \item The induced map 
    \begin{align*}
        \gr^\nu_F \Gamma_*: \gr^\nu_F CH^j(X)_\mathbb{Q} \rightarrow \gr^\nu_F CH^{j+\ell-\dim X}(Y)_\mathbb{Q}
    \end{align*}
    vanishes if the map of cohomology groups
    \begin{align*}
        [\Gamma]_*: H^{2j-\nu}(X, \mathbb{Q}) \rightarrow H^{2j-\nu + 2\ell -2 \dim X}(Y, \mathbb{Q})
    \end{align*}
    is zero on the Hodge pieces $H^{r,s}$ of $H^{2j-\nu}(X, \mathbb{Q})$ for $s \leq j - \nu$.
\end{enumerate}
\end{conjecture}

Conjecture \ref{conj bb} implies Conjecture \ref{intro-conj BB} by taking $\Gamma = \Delta_X$. Indeed, for every $1 \leq \nu \leq n$, condition $(4)$ implies 
\begin{align*}
    \gr_F^\nu CH_0(X)_{\mathbb Q}=0
\end{align*}
whenever $0 = H^{n,n-\nu}(X) \simeq H^0(X, \Omega_X^\nu)$. Thus, if $H^0(X,\Omega_X^\nu)=0$ for every $\nu > 0$, then all positive graded pieces of the filtration on $CH_0(X)_{\mathbb Q}$ vanish, in which case $CH_0(X)_{\mathbb Q}\simeq\mathbb Q$. Finally, Lemma~\ref{lem CH_0 Q} upgrades this rational isomorphism to $CH_0(X) \simeq \mathbb{Z}$.

We will need a stronger version of Conjecture~\ref{conj bb} to treat singular varieties. Consider a semisimplicial variety 
\[
\begin{tikzcd}[row sep=1.4em, column sep=3.5em]
X_\bullet = \bigg( \cdots X_2 \arrow[r] \arrow[r, shift left=0.9ex] \arrow[r, shift right=0.9ex] & 
X_1 \arrow[r] \arrow[r, shift left=0.6ex] &
X_0 \bigg),
\end{tikzcd}
\]
where each $X_a$ is smooth and projective. Usually, $X_\bullet$ will be a hyperresolution of a singular variety.

Taking the alternating sum of face maps gives rise to two natural complexes of Chow groups. Pullback along the face maps gives the cochain complex
\begin{align*}
    CH^m(X_0)_\mathbb{Q} \rightarrow CH^m(X_1)_\mathbb{Q} \rightarrow CH^m(X_2)_\mathbb{Q} \rightarrow \cdots,
\end{align*}
while proper pushforward gives the chain complex
\begin{align*}
    \ldots \rightarrow CH_m(X_2)_\mathbb{Q} \rightarrow CH_m(X_1)_\mathbb{Q} \rightarrow CH_m(X_0)_\mathbb{Q}.
\end{align*}

Conjecture~\ref{conj bb} gives a criterion for an individual morphism of Chow groups to vanish on a given graded piece of the Bloch--Beilinson filtration. We seek a corresponding criterion for either of the preceding complexes to be exact. This leads to the following conjecture. 

We recall that each $H_i(X_a, \mathbb{Q})$ carries a pure Hodge structure of weight $-i$, and it is dual to the pure Hodge structure on $H^i(X_a, \mathbb{Q})$; see \cite[Definition $1.17$]{peters steenbrink}.

\begin{conjecture}\label{conj bb exact}
Let $X_\bullet$ be a semisimplicial variety whose terms are smooth and projective. Fix integers $a, m \geq 0$.
\begin{enumerate}
    \item \emph{Cohomological Version.} The cochain complex
        \begin{align*}
            \cdots \rightarrow CH^m(X_{a-1})_\mathbb{Q} \rightarrow  CH^m(X_a)_\mathbb{Q} \rightarrow CH^m(X_{a+1})_\mathbb{Q} \rightarrow \cdots 
        \end{align*}
        is exact in the middle if the complex of pure Hodge structures
        \begin{align*}
            \cdots \rightarrow H^{2m-\nu}(X_{a-1}, \mathbb{C}) \rightarrow H^{2m-\nu}(X_{a}, \mathbb{C}) \rightarrow H^{2m-\nu}(X_{a+1}, \mathbb{C}) \rightarrow \cdots 
        \end{align*}
        is exact on every Hodge piece $H^{r,s} H^{2m-\nu}(X_a, \mathbb{C})$ of degree $r \geq m$ for all $\nu \geq 0$. Dually, this is equivalent to 
        \begin{align*}
            \cdots \rightarrow H_{2m-\nu}(X_{a+1}, \mathbb{C}) \rightarrow H_{2m-\nu}(X_{a}, \mathbb{C}) \rightarrow H_{2m-\nu}(X_{a-1}, \mathbb{C}) \rightarrow \cdots 
        \end{align*}
        being exact on every Hodge piece $H^{r,s} H_{2m-\nu}(X_a,\mathbb{C})$ of degree $r \leq -m$ for all $\nu \geq 0$.
        
    \item \emph{Homological Version.} The chain complex
        \begin{align*}
            \cdots \rightarrow CH_m(X_{a+1})_\mathbb{Q} \rightarrow  CH_m(X_a)_\mathbb{Q} \rightarrow CH_m(X_{a-1})_\mathbb{Q} \rightarrow \cdots 
        \end{align*}
        is exact in the middle if the complex of pure Hodge structures
        \begin{align*}
            \cdots \rightarrow H_{2m+\nu}(X_{a+1}, \mathbb{C}) \rightarrow H_{2m+\nu}(X_{a}, \mathbb{C}) \rightarrow H_{2m+\nu}(X_{a-1}, \mathbb{C}) \rightarrow \cdots 
        \end{align*}
        is exact on every Hodge piece $H^{r,s} H_{2m+\nu}(X_a, \mathbb{C})$ of degree $r \geq -m$ for all $\nu \geq 0$. Dually, this is equivalent to
        \begin{align*}
            \cdots \rightarrow H^{2m+\nu}(X_{a-1}, \mathbb{C}) \rightarrow H^{2m+\nu}(X_{a}, \mathbb{C}) \rightarrow H^{2m+\nu}(X_{a+1}, \mathbb{C}) \rightarrow \cdots 
        \end{align*}
        being exact on every Hodge piece $H^{r,s} H^{2m+\nu}(X_a, \mathbb{C})$ of degree $r \leq m$ for all $\nu \geq 0$.
\end{enumerate}
\end{conjecture}

The cohomological version was formulated in the more general setting of higher Chow groups in \cite[Conjecture~$19.5$]{bps}. Both versions are predicted by the same conjectural theory of mixed motives that underlies the classical Bloch--Beilinson conjecture, Conjecture~\ref{conj bb}.

Assuming Conjecture \ref{conj bb exact}, we obtain a precise relationship between zero-cycles and holomorphic forms for singular varieties. 

\begin{corollary}\label{cor CH_0 = Z}
Assume Conjecture \ref{conj bb exact}, and let $X$ be a connected projective variety.
\begin{enumerate}
    \item If $\mathbb{H}^0(X, \underline{\Omega}_X^p) = 0$ for every $p > 0$, then
    \begin{align*}
        CH_0(X) \simeq \mathbb{Z}.
    \end{align*}

    \item If $\mathbb{H}^p(X, \underline{\Omega}_X^0) = 0$ for every $p > 0$, then
    \begin{align*}
        RCH_0(X)_\mathbb{Q} \simeq \mathbb{Q}[0],
    \end{align*}
    concentrated in degree zero.
\end{enumerate}
\end{corollary}

This corollary motivats the singular Bloch-Beilinson Conjecture \ref{intro-conj singular bb} stated in the introduction. Indeed, the converse of $(1)$ is established in Theorem~\ref{intro-thm BB}. We expect the converse of $(2)$ to hold as well, although it remains open.

\begin{proof}
Take a hyperresolution $\epsilon_\bullet: X_\bullet \rightarrow X$. Denote by $\widetilde{X}_\bullet$ the same semisimplicial variety, but augmented with a point $\widetilde{X}_{-1} = \mathop{\mathrm{Spec}}\mathbb{C}$ at the end. 

$(1)$ Lemma \ref{lem RCH_0} gives an identification
\begin{align*}
    CH_0(X) = \ker \bigg( CH_0(X_0) \rightarrow CH_0(X_1) \bigg).
\end{align*}
Moreover, Lemma \ref{lem CH_0 Q} shows that $CH_0(X) \simeq \mathbb{Z}$ if $CH_0(X)_\mathbb{Q} \simeq \mathbb{Q}$.

Therefore, it suffices to prove that the augmented chain complex
\begin{align*}
    CH_0 \big( \widetilde{X}_{1} \big ) \rightarrow CH_0 \big( \widetilde{X}_{0} \big ) \rightarrow CH_0 \big( \widetilde{X}_{-1} \big )
\end{align*}
is exact in the middle. By Conjecture \ref{conj bb exact}$(2)$, this holds if the complex
\begin{align*}
    H^p \big ( \widetilde{X}_{-1}, \mathscr{O}_{\widetilde{X}_{-1}} \big) \rightarrow H^p \big ( \widetilde{X}_{0}, \mathscr{O}_{\widetilde{X}_{0}} \big) \rightarrow H^p \big ( \widetilde{X}_{1}, \mathscr{O}_{\widetilde{X}_{1}} \big) 
\end{align*}
is exact for all $p \geq 0$. But the strictness of the Hodge filtration together with the weight spectral sequence \eqref{eqn weights} shows that this middle cohomology is isomorphic to
\begin{align*}
    \gr_{0}^F \gr_p^W \widetilde{H}^p(X, \mathbb{C}) \simeq \gr_{-p}^F \widetilde{H}^p(X, \mathbb{C}).
\end{align*}
In degree $p > 0$, this is precisely the Hodge-Du Bois piece $\mathbb{H}^0(X, \underline{\Omega}_X^p)$ by \eqref{Hodge-dR ss}. This proves the claim.

$(2)$ We argue similarly. By Lemma \ref{lem CH_0 Q}, $RCH_0(X)_\mathbb{Q} \simeq \mathbb{Q}[0]$ if and only if the following complex is exact:
\begin{align*}
    \ldots \rightarrow CH_0 \big( \widetilde{X}_{1} \big )_\mathbb{Q} \rightarrow CH_0 \big( \widetilde{X}_{0} \big )_\mathbb{Q} \rightarrow CH_0 \big( \widetilde{X}_{-1} \big )_\mathbb{Q} 
\end{align*}
By Conjecture \ref{conj bb exact}$(2)$, this holds if the augmented Hodge complex
\begin{align*}
    H^p \big ( \widetilde{X}_{-1}, \mathscr{O}_{\widetilde{X}_{-1}} \big) \rightarrow H^p \big ( \widetilde{X}_{0}, \mathscr{O}_{\widetilde{X}_{0}} \big) \rightarrow H^p \big ( \widetilde{X}_{1}, \mathscr{O}_{\widetilde{X}_{1}} \big) \rightarrow \ldots
\end{align*}
is exact for all $p \geq 0$. By the argument above, the $\ell$-th cohomology of this complex is
\begin{align*}
    \gr_0^F \gr_p^W \widetilde{H}^{p + \ell}(X, \mathbb{C}).
\end{align*} 
Therefore, they all vanish if and only if $\gr_0^F H^p(X, \mathbb{C}) \simeq \mathbb{H}^0(X, \underline{\Omega}_X^p) = 0$ for all $p > 0$.
\end{proof}

\subsection{Chow-theoretic singularity conditions}\label{subsection Chow Dm}

In this section, we analyze the singularity conditions that arise naturally in the rigidity theorems of this paper. This leads to a hierarchy of Chow-theoretic singularity conditions, parallel to the Hodge-theoretic conditions $D_m$ for $m \geq 0$. 

Conditions $D_0$ and $E_0$ appear in Theorems \ref{thm Hodge general} and \ref{thm Hodge singular}, respectively. Both are Hodge-theoretic weakenings of rational singularities, and Proposition \ref{prop hiatt} shows that $D_0$ implies $E_0$. We begin by establishing pointwise characterizations of these conditions, which will facilitate comparison with their Chow-theoretic analogues.

\begin{proposition}\label{prop pointwise}
Let $X$ be a variety, and let $f: Y \rightarrow X$ be a (equivalently, any) resolution of singularities. Then:
\begin{enumerate}
    \item $X$ satisfies condition $D_0$ if and only if every fiber $F$ is connected and satisfies
    \begin{align*}
        \mathbb{H}^p(F, \underline{\Omega}_F^0) = 0 \qquad \text{ for every } p > 0.
    \end{align*}
    \item $X$ satisfies condition $E_0$ if and only if every fiber $F$ is connected and satisfies
    \begin{align*}
        \mathbb{H}^0(F, \underline{\Omega}_F^p) = 0 \qquad \text{ for every } p > 0.
    \end{align*}
\end{enumerate} 
\end{proposition}

\begin{proof}
Let $f: Y \rightarrow X$ be a resolution of singularities. 

$(1)$ By the characterization of condition $D_0$ in Proposition \ref{prop Dm}, it is enough to show that, for every closed point $x \in X$ and $q \geq 0$, there is an isomorphism
\begin{align}\label{eqn D0 iso}
     \gr_{-n}^F H^{q}_{\{ x \}} \big(  X, \mathbb{C} \big) \simeq \mathbb{H}^{2n-q}\big( f^{-1}(x), \underline{\Omega}_{f^{-1}(x)}^0)^\vee.
\end{align}
Indeed, for the only if direction, the equality $\mathbb{H}^0(F, \underline{\Omega}_F^0) = \mathbb{C}$ implies that every connected component of $X$ is irreducible, so the hypotheses of Proposition \ref{prop Dm} are satisfied. Conversely, the if direction also follows directly from Theorem \ref{intro-thm Omega^0}, although we do not need this stronger result here.

To prove the isomorphism \eqref{eqn D0 iso}, fix $x \in X$ and set $F = f^{-1}(x)$. Let $L_x$ denote the link of $x$ in $X$. There is a standard identification $H^{q}_{\{ x \}} \big(  X, \mathbb{C} \big) \simeq \widetilde{H}^{q-1} \big (L_x, \mathbb{C} \big)$. Moreover, if $N$ is a sufficiently small tubular neighborhood of $F$ in $Y$, \cite[Theorem $6.9$]{peters steenbrink} gives a long exact sequence for the pair $(N, F)$:
\begin{align}\label{eqn les link}
    \ldots \rightarrow H^{q-1}(F, \mathbb{C}) \rightarrow H^{q-1}(L_x, \mathbb{C}) \rightarrow H^{q}_F(Y, \mathbb{C}) \rightarrow H^q(F, \mathbb{C}) \rightarrow \ldots
\end{align}

Since the maps in this sequence are morphisms of mixed Hodge structures, strictness of the Hodge filtration allows us to pass to $\gr_{-n}^F$. Because $\dim F < n$, the outer terms of the resulting sequence vanish, so there is an isomorphism
\begin{align*}
    \gr_{-n}^F H^{q-1}(L_x, \mathbb{C}) \simeq \gr_{-n}^F H^q_F(Y, \mathbb{C}).
\end{align*}

Poincaré–Lefschetz duality identifies the dual of the latter group with $\gr_0^F H^{2n-q} (F, \mathbb{C})$. By degeneration of the Hodge-de Rham spectral sequence \eqref{Hodge-dR ss}, this is isomorphic to $\mathbb{H}^{2n-q}(F, \underline{\Omega}_F^0)$. Combining these identifications proves \eqref{eqn D0 iso}.

$(2)$ As in the proof of Theorem \ref{thm Hodge singular}, let $M_f^\bullet \in D^b( \mathrm{MHM}(X))$ be the cone in the distinguished triangle
\begin{align}\label{eqn defn Mf}
    \mathbb{Q}_X^H \rightarrow f_* \mathbb{Q}_Y^H \rightarrow M_f^\bullet \xrightarrow[]{+1}.
\end{align}

By Remark \ref{rmk E_0 char}, $X$ satisfies condition $E_0$ if and only if 
\begin{align*}
    \mathcal{H}^j \gr_{-p}^F \dr ( M_f^\bullet) = 0 \hspace{.5cm} \text{ whenever } j-p \leq 0.
\end{align*}
Propositions \ref{prop mhm 2} and \ref{prop mhm 2 converse} show that this is equivalent to requiring, for every closed point with embedding $i_x^*: \{ x \} \hookrightarrow X$,
\begin{align*}
    \mathcal{H}^j \gr_{-p}^F \dr ( i_x^* M_f^\bullet) = 0 \hspace{.5cm} \text{ whenever } j-p \leq 0.
\end{align*}
Applying $i_x^*$ to \eqref{eqn defn Mf} and using the base-change result of Proposition \ref{prop base change}, we obtain
\begin{align*}
    \mathcal{H}^j \gr_{-p}^F \dr ( i_x^* M_f^\bullet) \simeq \gr_{-p}^F \widetilde{H}^j(F, \mathbb{C}),
\end{align*}
where $F = f^{-1}(x)$. The right hand side vanishes automatically when $j < p$. When $j = p = 0$, it vanishes precisely when $F$ is connected. For $j = p > 0$, it is identified with the Hodge-de Rham piece $\mathbb{H}^0(F, \underline{\Omega}_F^p)$. This completes the proof.
\end{proof}

The pointwise characterization of condition $E_0$ above suggests that it is equivalent to the notion of $CH_0$-trivial singularities introduced in Lemma \ref{lem CH_0 trivial}.

\begin{corollary}\label{cor CH_0 vs E_0}
The following statements hold.
\begin{enumerate}
    \item If $X$ has $CH_0$-trivial singularities, then $X$ satisfies condition $E_0$.
    \item Assuming Conjecture \ref{conj bb exact}, the converse holds: if $X$ satisfies condition $E_0$, then it has $CH_0$-trivial singularities
\end{enumerate}
\end{corollary}

\begin{proof}
$(1)$ By Lemma \ref{lem CH_0 trivial}, $X$ has $CH_0$-trivial singularities if and only if, for a resolution of singularities $f: Y \rightarrow X$, every fiber $F$ satisfies
\begin{align*}
    CH_0(F) \simeq \mathbb{Z}.
\end{align*}
We can always choose the resolution so that $F$ is an snc divisor, hence is equidimensional. Theorem \ref{intro-thm BB} then gives
\begin{align*}
        \mathbb{H}^0(F, \underline{\Omega}_F^p) = 0 \qquad \text{ for every } p > 0.
\end{align*}
The result follows from Proposition \ref{prop pointwise}$(2)$. 

$(2)$ Assuming Conjecture \ref{conj bb exact}, the converse follows similarly by applying Corollary \ref{cor CH_0 = Z}$(1)$.
\end{proof}

Since rational singularities imply condition $E_0$, Corollary \ref{cor CH_0 vs E_0}$(2)$ suggests the following concrete conjecture.

\begin{conjecture}\label{conj ratl implies CH_0}
If $X$ has rational singularities, then $X$ has $CH_0$-trivial singularities. In particular, if $f: Y \rightarrow X$ is a resolution of singularities, then the induced map
\begin{align*}
    CH_0(Y) \rightarrow CH_0(X)
\end{align*}
is an isomorphism.
\end{conjecture}

For projective cones over smooth varieties, Conjecture~\ref{conj ratl implies CH_0} specializes to the generalized Bloch conjecture Conjecture~\ref{intro-conj BB}. Indeed, let $X$ be a projective cone over a smooth variety $Z$ satisfying $H^0(Z, \Omega_Z^p) = 0$ for all $p > 0$. Then $X$ has rational singularities whenever the conormal bundle is sufficiently positive by \cite[Proposition B.$3$]{ks}, and it is clear that $CH_0(X) \simeq \mathbb{Z}$. On the other hand, if $f: Y \rightarrow X$ is the blow-up of the vertex, then $CH_0(Y) \simeq CH_0(Z)$. 

We now list the known cases of Conjecture \ref{conj ratl implies CH_0}.

\begin{proposition}
Conjecture \ref{conj ratl implies CH_0} holds when $X$ is
\begin{enumerate}
    \item of klt type;
    \item of dimension at most two;
    \item a threefold, assuming Bloch's Conjecture for surfaces.
\end{enumerate}
\end{proposition}

\begin{proof}
$(1)$ follows from Lemma \ref{lem RCH_0 implications}, proved below.

$(2)$ If $\dim X < 2$, then rational singularities are smooth, so there is nothing to prove. We may therefore assume that $X$ is a surface. Let $f: Y \rightarrow X$ be a resolution of singularities. By Lemma \ref{lem CH_0 trivial}, it suffices to show that $CH_0(f^{-1}(x)) \simeq \mathbb{Z}$ for each fiber $F = f^{-1}(x)$. By Lemma \ref{lem wft}, we may assume that $F$ has snc support. Since $X$ has rational singularities, $F$ must be a connected tree of rational curves. Hence $CH_0(F) \simeq \mathbb{Z}$, as required.

$(3)$ As above, it suffices to prove that $CH_0(F) \simeq \mathbb{Z}$ for every snc fiber $F = \cup_i S_i$ of a resolution $f: Y \rightarrow X$. By Proposition \ref{prop pointwise}, we have vanishing
\begin{align*}
    0 = H^2(F, \mathscr{O}_F) \simeq \bigoplus_i H^2(S_i, \mathscr{O}_{S_i} ).
\end{align*}
Therefore, Bloch's Conjecture implies that $CH_0(S_i)_{\mathrm{hom}} \simeq \mathrm{Alb}(S_i)$ for every $i$.

There is a presentation
\begin{align*}
    CH_0(F) \simeq \coker \bigg( \bigoplus_{i,j} CH_0(C_{ij}) \xrightarrow{\alpha} \bigoplus_i CH_0(S_i) \bigg),
\end{align*}
where $C_{ij} = S_i \cap S_j$ denotes the pairwise intersection curve. Using the degree maps, $\alpha$ fits into a commutative diagram of the form
\[
\begin{tikzcd}[column sep=large, row sep=large]
\displaystyle\bigoplus_{i,j} \operatorname{Alb}(C_{ij})
  \arrow[r, "\beta"] \arrow[d]
&
\displaystyle\bigoplus_i \operatorname{Alb}(S_i)
  \arrow[d]
\\
\displaystyle\bigoplus_{i,j} \mathrm{CH}_0(C_{ij})
  \arrow[r, "\alpha"] \arrow[d, "\deg"]
&
\displaystyle\bigoplus_i \mathrm{CH}_0(S_i)
  \arrow[d, "\deg"]
\\
\displaystyle\bigoplus_{i,j} H_0(C_{ij},\mathbb{Z})
  \arrow[r, "\gamma"]
&
\displaystyle\bigoplus_i H_0(S_i,\mathbb{Z}),
\end{tikzcd}
\]
where the columns are short exact sequences. We have $\coker \gamma \simeq H_0( \mathcal{D}(F), \mathbb{Z})$ for the dual complex $\mathcal{D}(F)$ of $F$. Since $X$ is normal, the fiber $F$, and hence its dual complex $\mathcal{D}(F)$, is connected, so $\coker \gamma \simeq \mathbb{Z}$.

The snake lemma then shows that it suffices to prove $\beta$ is surjective. We view the Albanese as a quotient $\mathrm{Alb}(S_i) \simeq H^0(S_i, \Omega_{S_i}^1)^\vee / H_1(S_i, \mathbb{Z})$ and $\mathrm{Alb}(C_{ij}) \simeq H^0(C_{ij}, \Omega_{C_{ij}}^1)^\vee / H_1(C_{ij}, \mathbb{Z})$. It is therefore enough to prove that the induced map of vector spaces
\begin{align*}
    \bigoplus_{i,j} H^0(C_{ij}, \Omega_{C_{ij}}^1)^\vee \rightarrow \bigoplus_i H^0(S_i, \Omega_{S_i}^1)^\vee
\end{align*}
is surjective. Equivalently, its dual restriction map
\begin{align*}
    \bigoplus_{i} H^0(S_i, \Omega_{S_i}^1) \rightarrow \bigoplus_{i,j} H^0(C_{ij}, \Omega_{C_{ij}}^1)
\end{align*}
must be injective. 

By the degeneration of the weight spectral sequence \ref{eqn weights}, the kernel of this map is identified with 
\begin{align*}
    \gr_0^F \gr_1^W H^1(F, \mathbb{C}) \simeq \gr_1^W \mathbb{H}^1(F, \underline{\Omega}_F^0).
\end{align*}
This group vanishes by Proposition \ref{prop pointwise}, and hence the desired injectivity follows.
\end{proof}

We now introduce a Chow-theoretic analogue of condition $D_0$, which mirrors the relationship between condition $E_0$ and $CH_0$-trivial singularities.

\begin{definition}
Let $X$ be a variety, and let $f: Y \rightarrow X$ be a (equivalently, any) resolution of singularities. We say that $X$ has \emph{$RCH_0$-trivial singularities} if every fiber $F$ satisfies
\begin{align*}
    RCH_0(F)_\mathbb{Q} \simeq \mathbb{Q}[0],
\end{align*}
concentrated in degree zero.
\end{definition}

The following lemma records what is known about $RCH_0$-trivial singularities.

\begin{lemma}\label{lem RCH_0 implications}
The following statements hold.
\begin{enumerate}
    \item If $X$ is of klt type, then $X$ has $RCH_0$-trivial singularities.
    \item If $X$ has $RCH_0$-trivial singularities, then $X$ has $CH_0$-trivial singularities.
    \item Assuming Conjecture \ref{conj bb exact}, condition $D_0$ implies $RCH_0$-trivial singularities.
\end{enumerate}
\end{lemma}

\begin{proof}
Let $f: Y \rightarrow X$ be a resolution of singularities.

$(1)$ By \cite[Proposition $6.13$]{burke dual}, $RCH_0$ is an RC-invariant functor. Theorem G of \emph{loc. cit.}, applied to $f$, proves the claim.

$(2)$ Let $F$ be a fiber of $f$. By Lemma \ref{lem RCH_0}, we have 
\begin{align*}
    CH_0(F)_\mathbb{Q} \simeq H_0(RCH_0(F)_\mathbb{Q}) \simeq \mathbb{Q}.
\end{align*}
Lemma \ref{lem CH_0 Q} then implies that $CH_0(F) \simeq \mathbb{Z}$. Hence, $X$ has $CH_0$-trivial singularities.

$(3)$ Let $F$ be a fiber of $f$. By Proposition \ref{prop pointwise}$(1)$,
\begin{align*}
        \mathbb{H}^p(F, \underline{\Omega}_F^0) = 0 \qquad \text{ for every } p > 0.
\end{align*}
Assuming Conjecture \ref{conj bb exact}, this implies that $RCH_0(F)_\mathbb{Q} \simeq \mathbb{Q}[0]$ by Corollary \ref{cor CH_0 = Z}$(2)$.
\end{proof}

Although we cannot prove a direct comparison between $RCH_0$-trivial singularities and condition $D_0$, we expect these two notions to be equivalent. We also do not know whether $RCH_0$-trivial singularities implies that, for a resolution $f$ and closed subvariety $Z \subset X$, there is an isomorphism 
\begin{align*}
RCH_0\big( f^{-1}(Z) \big)_\mathbb{Q} \simeq RCH_0 \big( Z \big)_\mathbb{Q}.
\end{align*}

The following examples show that rational coefficients are essential in the definition of $RCH_0$-trivial singularities. The first shows that the integral analogue of Conjecture \ref{intro-conj CH_0} is false, while the second exhibits a rational singularity that is not $RCH_0$-trivial with integral coefficients.

\begin{example}[{\bf Degeneration of Enriques surfaces}]\label{exmp enriques}
Consider a type III degeneration of Enriques surface $f: Y \rightarrow C$ as in \cite{morrison}. The general fiber is an Enriques surface $S$, so $CH_0(S) \simeq \mathbb{Z}$ by \cite{bloch conj}. On the other hand, the special fiber $E$ is snc with rational strata and dual complex $\mathcal{D}(E) \simeq \mathbb{RP}^2$. One computes that
\begin{align*}
    H_1 \big( RCH_0(E) \big) \simeq H_1 \big( \mathbb{RP}^2, \mathbb{Z} \big) \simeq \mathbb{Z}/2 \mathbb{Z},
\end{align*}
so Conjecture \ref{intro-conj CH_0} fails integrally.
\end{example}

\begin{example}[{\bf Threefold with rational singularities}]\label{exmp RP^2}
By \cite[Example $8.3$]{payne}, there exists a threefold $X$ with an isolated rational singularity $x \in X$, and a log resolution $f: Y \rightarrow X$ whose exceptional divisor $E$ has rational strata and dual complex $\mathcal{D}(E) \simeq \mathbb{RP}^2$. Then
\begin{align*}
    H_1 \big( RCH_0(E) \big) \simeq H_1 \big( \mathbb{RP}^2, \mathbb{Z} \big) \simeq \mathbb{Z}/2 \mathbb{Z},
\end{align*}
so $X$ is not $RCH_0$-trivial in the stronger sense with integral coefficients.
\end{example}

Our goal for the remainder of this section is to formulate and justify a higher notion of $RCH_m$-trivial singularities analogous to condition $D_m$. This requires a pointwise characterization of condition $D_m$, which we give in Proposition \ref{prop Dm}. By contrast, for $m > 0$, we do not know of a pointwise characterization of condition $E_m$, and consequently we do not have a higher notion of $CH_m$-triviality.

We begin with some notation. Let $E$ be a simple normal crossings divisor in a smooth variety $Y$ of dimension $n$. For $k \geq 0$, denote by $E^{[k]}$ the disjoint union of the $(k+1)$-fold intersections of the irreducible components of $E$.

For each $j \geq 0$, consider the complex
\begin{align*}
    \kappa_j := \bigg( \ldots \rightarrow CH_j(E^{[1]})_\mathbb{Q} \rightarrow CH_j(E^{[0]})_\mathbb{Q} \rightarrow CH^{n-j}(E^{[0]})_\mathbb{Q} \rightarrow CH^{n-j}(E^{[1])})_\mathbb{Q} \rightarrow \ldots \bigg).
\end{align*}
We place $CH_j(E^{[0]})_\mathbb{Q}$ in degree zero. The differentials in negative degrees are given by the usual alternating sums of pushforward maps, while those in positive degrees are given by alternating sums of pullback maps.

If we denote by $\iota: E^{[0]} \hookrightarrow Y$ the natural inclusion, then the middle differential is the composition 
\begin{align*}
    CH_j(E^{[0]})_\mathbb{Q} \xrightarrow{\iota_*} CH_j(Y)_\mathbb{Q} \simeq CH^{n-j}(Y)_\mathbb{Q} \xrightarrow{\iota^*} CH^{n-j}(E^{[0]})_\mathbb{Q}.
\end{align*}
It is straightforward to show that the differential squares to zero, so $\kappa_j$ is indeed a complex.

\begin{definition}\label{defn RCH_m}
Let $X$ be a variety of pure dimension $n$. We say that $X$ has \emph{$RCH_m$-trivial singularities} if, for every closed point $x \in X$, there exists a resolution of singularities $f: Y \rightarrow X$ such that $E = f^{-1}(x)_{\mathrm{red}}$ is a simple normal crossings divisor, and 
\[
\kappa_j \text{ is null-homotopic for every }0<j \leq m,
\qquad
\kappa_0\simeq \mathbb Q[0].
\]
\end{definition}

When $m = 0$, the complex $\kappa_0$ is identified with $RCH_0(E)_\mathbb{Q}$. Thus, this definition recovers the notion of $RCH_0$-triviality introduced above. More generally, we expect $RCH_m$-trivial singularities to be equivalent to condition $D_m$.

\begin{conjecture}\label{conj RCHm Dm}
A variety of pure dimension has $RCH_m$-trivial singularities if and only if it satisfies condition $D_m$.
\end{conjecture}

The following example shows that the forward implications holds for cones over a smooth projective base.

\begin{example}[{\bf Projective cone}]\label{exmp cone}
Let $Z$ be a smooth projective variety of dimension $n-1$, and let $X$ be a projective cone over $Z$. Then
\[
X\text{ has }RCH_m\text{-trivial singularities}
\quad\Longleftrightarrow\quad
CH_j(Z)_{\mathbb Q}\simeq \mathbb Q
\quad\text{for every }0\le j\le m.
\]

On the other hand, for example by Proposition \ref{prop Dm} below,
\[
X\text{ satisfies condition }D_m
\quad\Longleftrightarrow\quad
h^{p,q}(Z)=h^{p,q}(\mathbb P^{n-1})
\quad\text{for all }0\le p\le m\text{ and }q\ge 0.
\]
Therefore, in this situation, \cite[Theorem $3.20$]{voisin ii} proves that $RCH_m$-trivial singularities implies condition $D_m$.
\end{example}

To investigate the opposite implication in Conjecture \ref{conj RCHm Dm}, we establish a pointwise characterization of condition $D_m$ in Proposition \ref{prop Dm} and compare it with the motivic heuristics suggested by Conjecture \ref{conj bb exact}.

We fix the following notation. Let $E$ be a simple normal crossings divisor in a smooth variety $Y$ of dimension $n$. For each $j \geq 0$, let $\Theta_j$ denote the following complex of pure Hodge structures of weight $j$:
\begin{align*}
     \ldots \rightarrow H_{2n-j} \big( E^{[1]}, \mathbb{C} \big)(-n) \rightarrow H_{2n-j} \big( E^{[0]}, \mathbb{C} \big)(-n) \rightarrow H^{j} \big( E^{[0]}, \mathbb{C} \big) \rightarrow H^{j} \big( E^{[1]}, \mathbb{C} \big) \rightarrow \ldots 
\end{align*}

As before, we place $H_{2n-j}(E^{[0]}, \mathbb{C})(-n)$ in degree zero. The differentials in negative degrees are the alternating sums of the pushforward maps, while those in positive degrees are the alternating sums of the pullback maps. The middle differential is the composition
\begin{align*}
    H_{2n-j} \big( E^{[0]}, \mathbb{C} \big)(-n) \xrightarrow[]{\iota_*} H_{2n-j}(Y, \mathbb{C} \big) (-n) \xrightarrow{\simeq} H^j(Y, \mathbb{C}) \xrightarrow{\iota^*} H^j(E^{[0]}, \mathbb{C}),
\end{align*}
where the middle isomorphism is given by Poincaré duality. 

For $j = p+q$, define $\Theta_{p,q} = \gr_{-p}^F \Theta_j$. A direct computation shows that the complex $\Theta_{p,q}$ has the form
\begin{align*}
     \ldots \rightarrow H^{q-2} \big( E^{[1]}, \Omega_{E^{[1]}}^{p-2} \big) \rightarrow H^{q-1} \big( E^{[0]}, \Omega_{E^{[0]}}^{p-1} \big) \rightarrow H^q \big( E^{[0]}, \Omega_{E^{[0]}}^p \big) \rightarrow H^q \big( E^{[1]}, \Omega_{E^{[1]}}^p \big) \rightarrow \ldots
\end{align*}

We now state our pointwise characterization of condition $D_m$.

\begin{proposition}\label{prop Dm}
Let $X$ be a variety of pure dimension $n$. Then $X$ satisfies condition $D_m$ if and only if, for every closed point $x \in X$, there exists a resolution of singularities $f: Y \rightarrow X$ such that with $E = f^{-1}(x)_{\mathrm{red}}$ a simple normal crossings divisor, and
\begin{align*}
    \Theta_{p,q}\text{ is null-homotopic for every }0\le p\le m,\ q\ge0,
\end{align*}
except that $\Theta_{0,0} \simeq \mathbb{C}[0]$.
\end{proposition}

\begin{proof}
Fix a closed point $x \in X$, and choose a resolution $f: Y \rightarrow X$ such that $E = f^{-1}(x)_{\mathrm{red}}$ is a simple normal crossings divisor. 

First consider the case $m = 0$. By the weight spectral sequence \eqref{eqn weights}, we have $\Theta_{0,0} \simeq \mathbb{C}[0]$ and $\Theta_{0,q}$ is null-homotopic for every $q > 0$ if and only if $E$ is connected and $\mathbb{H}^q(E, \underline{\Omega}_E^0) = 0$ for every $q > 0$. Thus, in degree zero, the proposition is a consequence of Proposition \ref{prop pointwise}$(1)$.

The argument in higher degrees is analogous. By the characterization of condition $D_m$ in Theorem \ref{thm Dm}, it suffices to show that, for fixed $p \geq 1$, 
\begin{align}\label{eqn theta equiv}
\left(
\operatorname{gr}_{p-n}^F H^j_{\{x\}}(X,\mathbb C)=0
\ \text{for all }j\ge 0
\right)
\Longleftrightarrow
\left(
\Theta_{p,q}\text{ is null-homotopic}
\ \text{for all }q\ge 0
\right).
\end{align}

Consider the commutative diagram
\begin{equation}\label{eqn comm}
\begin{tikzcd}
H_E^{j}(Y,\mathbb C)
    \arrow[r, "\alpha"]
    \arrow[d, "\gamma"', "\simeq"]
&
H^{j}(Y,\mathbb C)
    \arrow[r, "\beta"]
    \arrow[d, "\simeq"]
&
H^{j}(E,\mathbb C)
\\
H_{2n-j}(E,\mathbb C)(-n)
    \arrow[r]
&
H_{2n-j}^{BM}(Y,\mathbb C)(-n).
\end{tikzcd}
\end{equation}
The vertical arrows are induced by cap product with the fundamental class $ \cap [Y]$. They are isomorphisms by Poincaré--Lefschetz duality.

As in the proof of Proposition \ref{prop pointwise}, the long exact sequence \eqref{eqn les link} shows that the vanishing condition on the left-hand side of \eqref{eqn theta equiv} is equivalent to the requirement that
\begin{align*}
    \gr_{p-n}^F ( \beta \circ \alpha ) \text{ is an isomorphism for all } j \geq 0.
\end{align*}
Since $\gamma$ is an isomorphism, this is equivalent to saying
\begin{align*}
    \gr_{p-n}^F ( \beta \circ \alpha \circ \gamma^{-1} ) \text{ is an isomorphism for all } j \geq 0.
\end{align*}
After dualizing and replacing $j$ by $2n-j$, we obtain the equivalent condition
\begin{align*}
    \gr_{-p}^F (\beta \circ \alpha \circ \gamma^{-1}) \hspace{.2cm} \text{ is an isomorphism for all } j \geq 0.
\end{align*}

Recall that $H^j(E, \mathbb{C})$ has weights at most $j$, whereas $H_{2n-j}(E, \mathbb{C})(-n)$ has weights at least $(j-2n)+2n=j$. Thus, the morphism $\beta \circ \alpha \circ \gamma^{-1}$ can be an isomorphism only if the Hodge pieces $\gr_{-p}^F H^j(E, \mathbb{C})$ and $\gr_{n-p}^F H_{2n-j}(E, \mathbb{C})$ are pure.

Under this assumption, the weight spectral sequence \eqref{eqn weights} shows that applying $\gr_{-p}^F$ to 
\begin{align*}
    H^j \big( E^{[0]}, \mathbb{Q} \big) \rightarrow H^j \big( E^{[1]}, \mathbb{Q} \big) \rightarrow H^j \big(E^{[2]}, \mathbb{Q} \big) \rightarrow \ldots
\end{align*}
produces a complex quasi-isomorphic to $\gr_{-p}^F H^j(E, \mathbb{C})$. Similarly, applying $\gr_{n-p}^F$ to
\begin{align*}
    \ldots \rightarrow H_{2n-j} \big( E^{[2]}, \mathbb{C} \big) \rightarrow H_{2n-j} \big( E^{[1]}, \mathbb{C} \big) \rightarrow H_{2n-j} \big( E^{[0]}, \mathbb{C} \big)
\end{align*}
yields a complex quasi-isomorphic to $\gr_{n-p}^F H_{2n-j}(E, \mathbb{C})$.

Splicing these two complexes together along $\beta \circ \alpha \circ \gamma^{-1}$, we conclude that $\gr_{-p}^F ( \beta \circ \alpha \circ \gamma^{-1})$ is an isomorphism if and only if $\Theta_{p,j-p} = \gr_{-p}^F \Theta_j$ is null-homotopic. Since this holds for every $j \geq 0$, the proof is complete.
\end{proof}

We conclude by comparing the complexes $\Theta_{p,q}$ and $\kappa_j$ via Conjecture \ref{conj bb exact}. Both complexes are obtained by splicing a cohomological complex and a homological complex along a central morphism.

The dual formulation of Conjecture \ref{conj bb exact}$(2)$ predicts that, if the right side of $\Theta_{p,q}$,
\begin{align*}
    \gr_{-p}^F H^{p+q} \big (E^{[0]}, \mathbb{Q} \big) \rightarrow \gr_{-p}^F H^{p+q} \big (E^{[1]}, \mathbb{Q} \big) \rightarrow \gr_{-p}^F H^{p+q} \big (E^{[2]}, \mathbb{Q} \big) \rightarrow \ldots
\end{align*}
is null-homotopic for all $0 \leq p \leq m$ and $q \geq 0$, then the left side of $\kappa_j$,
\begin{align}\label{eqn Chow 1}
   CH_j(E^{[2]})_\mathbb{Q} \rightarrow CH_j(E^{[1]})_\mathbb{Q} \rightarrow CH_j(E^{[0]})_\mathbb{Q} 
\end{align}
is null-homotopic for all $0 \leq j \leq m$.

Similarly, the dual formulation of Conjecture \ref{conj bb exact}$(1)$ predicts that, if the left side of $\Theta_{p,q}$,
\begin{align*}
    \ldots \rightarrow \gr_{n-p}^F H_{2n-p-q} \big (E^{[2]}, \mathbb{Q} \big) \rightarrow \gr_{n-p}^F H_{2n-p-q} \big (E^{[1]}, \mathbb{Q} \big) \rightarrow \gr_{n-p}^F H_{2n-p-q} \big (E^{[0]}, \mathbb{Q} \big)
\end{align*}
is null-homotopic for all $0 \leq p \leq m$ and $q \geq 0$, then the right side $\kappa_j$,
\begin{align}\label{eqn Chow 2}
    CH^{n-j}(E^{[0]})_\mathbb{Q} \rightarrow CH^{n-j}(E^{[1])})_\mathbb{Q} \rightarrow CH^{n-j}(E^{[2]})_\mathbb{Q} \rightarrow \ldots
\end{align}
is null-homotopic for all $0 \leq j \leq m$.

This discussion does not address the central differentials. Nevertheless, it suggests that the triviality of $\Theta_{p,q}$ for all $0 \leq p \leq m$ and $q \geq 0$ should imply the triviality of $\kappa_j$ for all $j \leq m$. Thus, we expect condition $D_m$ to imply $RCH_m$-trivial singularities as in Conjecture \ref{conj RCHm Dm}.

\end{document}